\documentclass[12pt]{amsart}

\numberwithin{equation}{section}
\theoremstyle{plain}
\newtheorem{theorem}[equation]{Theorem}

\newtheorem{proposition}[equation]{Proposition}
\newtheorem{lemma}[equation]{Lemma}
\newtheorem{corollary}[equation]{Corollary}
\newtheorem{conjecture}[equation]{Conjecture}

\theoremstyle{remark}
\newtheorem{remark}[equation]{Remark}

\theoremstyle{definition}
\newtheorem{definition}[equation]{Definition}

\newcommand{\lra}{\longrightarrow}
\newcommand{\ra}{\rightarrow}
\newcommand{\restr}{\mbox{\Large \(|\)\normalsize}}

\newcommand{\B}{{\mathcal B}}

\newcommand{\D}{{\mathcal D}}

\newcommand{\G}{{\mathcal G}}

\renewcommand{\H}{\mathbb H}

\renewcommand{\L}{{\mathcal L}}

\newcommand{\R}{\mathbb R}

\newcommand{\T}{{\mathcal T}}

\newcommand{\bad}{\operatorname{Bad}}
\newcommand{\bv}{\operatorname{BV}}
\newcommand{\cent}{\operatorname{Center}}

\newcommand{\const}{\operatorname{const}}
\newcommand{\cut}{\operatorname{Cut}}

\newcommand{\fp}{\operatorname{FP}}
\newcommand{\gl}{\operatorname{GL}}
\newcommand{\good}{\operatorname{Good}}
\newcommand{\halfspace}{\operatorname{HS}}
\newcommand{\hs}{\operatorname{HS}}

\newcommand{\id}{\operatorname{id}}

\newcommand{\Lip}{\operatorname{Lip}}
\newcommand{\LIP}{\operatorname{LIP}}

\newcommand{\loc}{\operatorname{loc}}

\newcommand{\mass}{\operatorname{Mass}}

\newcommand{\meas}{\operatorname{\mathcal D}}

\newcommand{\on}{\:\mbox{\rule{0.1ex}{1.2ex}\rule{1.9ex}{0.1ex}}\:}

\newcommand{\per}{\operatorname{Per}}
\newcommand{\Per}{\operatorname{PER}}
\newcommand{\radon}{\operatorname{Radon}}

\newcommand{\super}{\operatorname{Slice}}
\newcommand{\supp}{\operatorname{Supp}}
\newcommand{\taut}{\operatorname{Taut}}

\newcommand{\Var}{\operatorname{VAR}}
\newcommand{\var}{\operatorname{Var}}

\def\D{\partial}
\newcommand{\al}{\alpha}
\def\de{\delta}
\def\De{\Delta}

\def\eps{\epsilon}
\def\ga{\gamma}

\def\la{\lambda}
\def\La{\Lambda}

\def\lra{\longrightarrow}

\def\om{\omega}

\def\ra{\rightarrow}

\def\si{\sigma}
\def\Si{\Sigma}

\def\be{\beta}

\def\defeq{:=}

\newcommand{\ol}{\overline}
\newcommand{\no}{\noindent}

\def\XXint#1#2#3{{\setbox0=\hbox{$#1{#2#3}{\int}$}
      \vcenter{\hbox{$#2#3$}}\kern-.5\wd0}}

\begin{document}

\begin{abstract}
This is one of a series of papers
examining the interplay between differentiation theory 
for Lipschitz maps, $X\ra V$, and bi-Lipschitz nonembeddability, 
where $X$ is a metric measure space and $V$ is a Banach
space. Here, we consider the case $V=L^1$,
where differentiability fails.
 We establish another kind of differentiability
for certain
$X$, including
${\R}^n$ and $\H$, the Heisenberg group
 with its Carnot-Caratheodory metric. 
It follows
that $\H$ does not bi-Lipschitz embed into $L^1$, 
as conjectured by J. Lee and A. Naor. 
When combined with their work,
this provides a natural counterexample to the Goemans-Linial conjecture
in theoretical computer science; the first such counterexample
was found by Khot-Vishnoi \cite{khotvishnoi}. 
A key ingredient in the proof of our main theorem is a new connection between
Lipschitz maps to $L^1$ and functions of bounded variation, 
which permits us to exploit recent work on the structure of
$\bv$ functions on the Heisenberg group \cite{italians1}. 
\end{abstract}

\title[Differentiating maps into $L^1$]{Differentiating maps into $L^1$, 
 and the geometry of $\bv$ functions}
\date{\today}
\author{Jeff Cheeger}
\address{J.C.\,: Courant Institute of Mathematical Sciences\\
       251 Mercer Street\\
       New York NY 10012}

\author{Bruce Kleiner}
\address{B.K.\,: Mathematics Department\\
         Yale University\\
             New Haven, CT 06520}
\thanks{The first author was partially supported by NSF Grant
DMS 0105128 and the second by NSF Grant DMS 0505610}
\maketitle

{\small\tableofcontents}

\section{Introduction}
\vskip4mm

\no
\subsection{Overview}

\mbox{}

\vspace{2mm}
The interplay between differentiability and bi-Lipschitz nonembeddability 
is the common theme of this paper and the papers \cite{sepdual,laaksoembed},
see also \cite{crannouncement}.
Specifically, we are concerned with the case in which the 
target is an infinite dimensional Banach space and the domain is
a complete metric measure space for which
the measure satisfies a doubling condition and a Poincar\'e inequality holds
in the sense of upper gradients; see \cite{HeKo}.  
Such metric measure spaces will be referred to as {\it PI spaces}.

In \cite{sepdual}, the differentiability theorem for
real valued Lipschitz functions on all PI spaces
and the resulting bi-Lipschitz nonembedding
theorem for certain PI spaces, proved in \cite{cheeger}
for finite dimensional Banach space targets, are extended to 
a class of infinite dimensional Banach space targets which admit
what we call {\it good finite dimensional approximation}.
Included in this class are
separable dual spaces.
Domains covered by this
nonembedding theorem include Bourdon-Pajot spaces,  Laakso spaces
(which are PI 
spaces of topological dimension $1$, for which the Hausdorff dimension 
can be any real number $>1$, \cite{laakso}) and the 
Heisenberg group, $\H$, with its Carnot-Caratheodory metric $d^\H$.

In the present paper, we examine maps, $X\to L^1
=L^1(\R,\L)$, where $\L$ denotes
Lebesgue measure on the real numbers $\R$.
We show that for a special class of PI spaces, including
$\R^k$ and $(\H,d^\H,\L)$,
despite the failure of the usual form of  differentiability for
Lipschitz maps to $L^1$,
a novel form of differentiability does in fact hold. As a direct
consequence, it follows that $(\H,d^\H)$
does not bi-Lipschitz embed in $L^1$. This proves
a conjecture of J. Lee and A. Naor, which provided the motivation for our work.
The significance of this conjecture
in the context of  theoretical computer science is briefly indicated below.

Differentiability is well known 
to fail even for Lipschitz maps $f:\R\to L^1$.
In particular, the differentiation 
theorem of \cite{sepdual} does not apply. Moreover,
 the nonembedding corollary turns
out not to hold.
In \cite{laaksoembed}, it is shown that members of a class of PI
spaces, including Laakso spaces,
 which satisfy the assumptions of the nonembedding theorem
of \cite{sepdual}, 
do embed in $L^1$.

\no
\subsection{ Rademacher's theorem and its descendents.}

\mbox{}

\vspace{2mm}Rademacher's  differentiation theorem 
states that a Lipschitz 
map, $f:\R^k\ra \R^l$, is differentiable almost everywhere. 
Hence, the geometry
of such Lipschitz maps becomes rigid at small scale. 
Specifically,
in the limit under suitable rescaling, the map becomes linear.
 The literature contains numerous extensions of this
result, in which either the domain, the target, or the class 
of maps is
generalized. Classical examples 
include almost everywhere approximate differentiability of 
Sobolev functions 
\cite{ziemer},
analysis of and on rectifiable sets, almost everywhere
differentiability of quasiconformal homeomorphisms between 
domains in $\R^k$,
and almost everywhere differentiability of Lipschitz maps 
from domains in $\R^k$, 
to certain Banach spaces which are said to possess the 
Radon-Nikodym
property; \cite[Chapter 5]{benlin}.  

In many of the recent results in this vein, a significant 
part of the achievement is 
to make sense of differentiation in a context where some 
component of the classical
setting is absent, e.g. the infinitesimal affine structure 
on the domain or target, or a good measure on the domain:

\noindent
$\bullet$ Pansu's differentiation theorem \cite{pansu} for
Lipschitz maps between graded nilpotent Lie groups.  Here,
 one cannot use
Euclidean rescaling; the
rescaling procedure has to be adapted to the grading on the 
groups.

\noindent
$\bullet$ Differentiation theory for real-valued
Lipschitz functions on Banach spaces with a separable dual 
\cite{linpreisssurvey}.
This requires replacing the classical notion of 
``almost everywhere'' by
something else.

\noindent
$\bullet$  Metric differentiation \cite{kirchheim,pauls}. 
 The target is an arbitrary metric space with no linear structure,
and differentiability is reformulated as a property 
of the pullback of the distance function from the target.

\noindent
$\bullet$  The differentiation theory developed in \cite{cheeger} 
for Lipschitz functions on PI spaces.
Typically, these carry no infinitesimal affine structure.
%\pagebreak

\noindent
\subsection{ Differentiation and bi-Lipschitz nonembeddability.}

\mbox{}

\vspace{2mm}Since differentiation theorems assert that the small-scale 
structure of maps
is very restricted, one can use them to show that certain  
mapping problems
have no solution.  For instance, it was observed by Semmes, 
\cite{semmes}, that
Pansu's differentiation theorem implies that
a Lipschitz map, $f:U\ra \R^k$,
where $U$ is an open subset of $(\H,d^\H)$,  
cannot be bi-Lipschitz.
Similarly, the differentiation theory of \cite{cheeger} 
was applied to give a unified proof
of  bi-Lipschitz nonembeddability in $\R^k$ of several 
families of spaces 
(\cite[Sect. 14]{cheeger}) including Carnot groups such as $(\H,d^\H$), 
 Laakso 
spaces, and the Bourdon-Pajot spaces (of \cite{bourpaj}).

The paper \cite{sepdual} extends  the differentiation 
theory of \cite{cheeger}
to Lipschitz maps $X\ra V$, where $(X,d^X,\mu)$ is a PI 
space and
$V$ is a Banach space with good finite dimensional
approximation, e.g. a separable Banach
space which is isomorphic to the dual of some Banach space.
As a consequence, the statement and proof of the  bi-Lipschitz 
nonembedding theorem 
of \cite[Sect. 14]{cheeger}, extend verbatim to separable
 dual targets.
In particular, this covers Lipschitz maps into arbitrary 
reflexive
spaces (separable or not), such as the $L^p$ space  
of an arbitrary measure space, $(Y,\nu)$,
for $1<p<\infty$, and also maps into the space,
$\ell^1$, of absolutely summable sequences.

In light of the above,  
the existence of bi-Lipschitz embeddings for certain domains 
and targets, implies the
nonexistence of a differentiation theorem. For example,
every metric space, $X$, admits a canonical isometric  
embedding into the space, $L^{\infty}(X)$,
the Kuratowski embedding, 
which assigns to $x\in X$,
the  function, $ d^X(x,\,\cdot\,)-d^X(x_0,\,\cdot\,)$, where $x_0\in X$
is a basepoint. Hence,
there cannot be {\it any} relevant differentiation theorem
 for maps into
$L^{\infty}$. In particular, the
procedure employed in the present paper to circumvent the failure
of the standard differentiation theorem for $L^1$ targets 
is useless when the target is $L^\infty$; see Remark \ref{kembed}.
\vskip2mm

\no
\subsection{ Failure of differentiability for Lipschitz maps to $\mathbf{L^1}$.}

\mbox{}

\vspace{2mm} For the target, $L^1(\R)$,  the failure of
differentiation theory is well-known, and is 
illustrated by the ``moving characteristic
function'' (cf. \cite{aronszajn}):
$$
f:[0,1]\ra L^1(\R)\, ,\quad f(t)\defeq \chi_{[0,t]\, },
$$
where $\chi_A$ denotes the characteristic function of 
a subset $A\subset\R$. Note that $f$ is actually
an isometric embedding. For this map, the difference
quotients at $t\in\R$ do not converge in $L^1$, but 
rather, when regarded
as measures, convergence weakly to the delta function, 
$\delta_t\not\in L^1$, concentrated at $t$.

\subsection{Bi-Lipschitz embedding of Laakso spaces in $L^1$.}

\mbox{}

\vspace{2mm}The main result of \cite{laaksoembed} states that members
of a class of spaces,
which includes the PI spaces of Laakso (as well as other 
interesting PI spaces)
admit bi-Lipschitz embeddings into $L^1$. According to  
\cite{sepdual}, these spaces do not bi-Lipschitz embed in
$\ell^1$. In particular,
for these domains,  the differentiation
result of \cite{sepdual} cannot be extended to $L^1$ 
in any form that is relevant to bi-Lipschitz nonembeddability.

To our knowledge, these 
spaces 
are the first examples of doubling metric spaces 
which bi-Lipschitz 
embed in $L^1$, but do not bi-Lipschitz embed in $\ell^1$.

\subsection{Bi-Lipschitz nonembedding in $L^1$; the Heisenberg group.}

\mbox{}

\vspace{2mm}The Heisenberg group, $(\H,d^H,\L)$, 
where $\L$ denotes Lebesgue measure,
is a PI space.
The motivation for  \cite{laaksoembed} and for the
present paper
came from the following conjecture of J. Lee and
A. Naor which is proved here.
\begin{conjecture}[Lee-Naor]
\label{quesassaf}
$(\H,d^\H)$, the Heisenberg group equipped with its 
Carnot-Caratheodory metric, does not
admit a bi-Lipschitz embedding into $L^1$.
\end{conjecture}

This conjecture arose from 
 \cite{leenaor}, in which it is shown that the nonexistence of such an
embedding would provide a natural counter-example to
 the Goemans-Linial conjecture of  theoretical
computer science;   for the first such counterexample, see \cite{khotvishnoi}. 
Very roughly, the point 
is that in some instances, questions in algorithm design, such as
the sparsest cut problem, could be solved if it were possible to
embed a certain class of finite metric spaces (those  with metrics of 
negative type) into $\ell^1$ with universally bounded 
bi-Lipschitz distortion, i.e. distortion independent of the particular
metric and the cardinality.

We now state a simplified version of our  
differentiation theorem. Let $e\in \H$ denote the identity
element.

\begin{theorem}
[Center collapse]
\label{thmmainfirst}
If $U\subset\H$ is an open subset, and $f:U\ra L^1$ is a Lipschitz map,
then for almost every point $x\in \H$, the map collapses in the
direction of the center of $\H$, i.e. 
\begin{equation}
\lim_{g\to e}\,\frac{\|f(gx)-f(x)\|_{L^1}}{d^\H(gx,x)}= 0\, ,\quad
\quad g\in\cent(\H)\, .
\end{equation}
\end{theorem}

Theorem \ref{thmmainfirst} implies that $f$ cannot 
be a bi-Lipschitz embedding, thus
proving Conjecture \ref{quesassaf}.  In particular:
\begin{corollary}
\label{cordoublingnotinellone}
There is a compact doubling metric space which does not bi-Lipschitz embed
in $L^1$.
\end{corollary}

To our knowledge, the Heisenberg group provides the first example of a
metric space with property stated in  Corollary \ref{cordoublingnotinellone}.

Two metric spaces, $W_1,\;W_2$, are called {\it quasi-isometric} if 
for some $D<\infty$, there
exist $D$-dense subsets, $\La_i\subset W_i$, and 
a bi-Lipschitz homeomorphism, $\La_1\ra \La_2$.

If $(W,d^W)$ is quasi-isometric to $(\H,d^\H)$, then the rescaled
sequence, $(W,i^{-1}d^W)$, converges in the 
(pointed) Gromov-Haudorff sense to a metric space bi-Lipschitz
homeomorphic to $(\H,d^\H)$ .

\begin{corollary}
\label{corqicor}
A metric space, $W$, that is quasi-isometric to $\H$
does not admit a bi-Lipschitz embedding in $L^1$.
\end{corollary}

Corollary \ref{corqicor} follows from Theorem \ref{thmmainfirst} by 
applying a general limiting argument \cite{heinmank,benlin}.
If the statement were
false, then by the theory of ultralimits, $(\H,d^\H)$  
would bi-Lipschitz embed in some Banach space,
$V$, which is an ultralimit of $L^1$ spaces. Then from  Kakutani's  
abstract characterization of $L^1$ spaces, \cite{kakutani}, it follows
that $V$ is itself an $L^1$ space; this contradicts Theorem \ref{thmmainfirst}.

The canonical example of  a space, $W$, to which 
Corollary \ref{corqicor}
applies is a Cayley graph $W$ for the integer
Heisenberg group, i.e. the subgroup of $\H$ for which
$a,b,c$ of (\ref{e:1.1}) below are integers. 
Recall that a   Cayley
graph $W$ for a group, $G$, is 
obtained by choosing a finite generating 
set $\Si\subset G$ and declaring
that two elements $g,g'\in G$ span an edge 
in $W$ if and only if
$g=g'\si$ for some $\si\in\Si$.  We equip $W$
with a $G$-invariant path metric.

Let $W$ denote some Cayley graph for the integer Heisenberg group, 
and for $k\geq 0$, let 
$$
W_k \defeq B_k(e) \subset W
$$
denote the combinatorial $k$-ball in $W$.

\begin{corollary}
\label{ul}
The sequence, $\{W_k\}$, is
a sequence of uniformly doubling finite graphs 
with uniformly bounded valence, 
which do not admit  embeddings into $L^1$ with uniformly bounded
bi-Lipschitz distortion.
\end{corollary}
\vskip2mm

\no
\subsection{Indication of proof.}
\label{subsecindication}

\mbox{}

\vspace{2mm}Our approach to differentiating maps to 
$L^1$ begins with
the equivalence between (pseudo-)metrics 
$d_f$,  on a set, $X$, which are induced by a map
$f:X\to L^1$, and
metrics which are representable as so-called ``cut metrics''.

 For now, the term {\it cut} just means  {\it subset}.
A cut, $E\subset X$, defines an
``elementary cut metric'', $d_E$, for
which $x_1,x_2$ have distance $1$ if either both points
lie in $E$ or neither point lies in $E$, and distance
$0$ otherwise.  A {\it cut metric}, $d_\Si$, is a superposition 
 of elementary cut metrics, with respect to a measure, $\Si$,
on the power set $2^X$;
$$
d_\Si(x_1,x_2)=\int_{2^X}d_E(x_1,x_2)\, d\Si\, .
$$
The measure, $\Si$, is called
a {\it cut measure}.

The basic fact (see e.g. Lemma 4.2.5 of \cite{dezalaur})
is that {\it any} metric, $d_f$, induced by a map, $f$,
from $X$ to an $L^1$ space,  can 
be realized as a cut
metric,  $d_{\Si_f}$, relative to a  cut measure, $\Si_f$,
canonically associated to $f$;
\begin{equation}
\label{cutmetrep}
d_f(x_1,x_2)=\int_{2^X}d_E(x_1,x_2)\, d\Si_f\, .
\end{equation}

In actuality, the set theoretic framework just described 
is not adequate for our subsequent purposes and
we will require a variant  in which
$X$ carries a $\sigma$-finite measure, $\mu$; see 
Section \ref{secl1mappings}.
However, 
 for the remainder of this subsection we will ignore this point.

Our main new  observation is that if $X$ is a PI space
and the map $f$  is Lipschitz, or more generally 
of bounded variation, then the cut measure 
$\Si$  will be supported on a very special subset
of $2^X$, namely on those, $E\subset 2^X$, with {\it finite perimeter};
see Section \ref{prelims} for the definition and some basic properties
of sets of finite perimeter.

Let $U\subset \H$ be  open  and let
$E\subset\H$ have finite perimeter in $U$.
Let $\per(E,U)\subset \rm{Radon}(U)$ denote the perimeter  measure 
of $E$ in $U$. By a recent structure theorem
in geometric measure theory,
for $\per(E,U)$ almost 
every point $x\in U$, when
one blows up $E$ at $x$, the resulting sequence of
characteristic functions
 converges in $L^1_{\loc}$
to a half-space; 
\cite{luigi1,luigi2,italians1,italians2}.  
Here a  {\it half-space} in the 
Heisenberg group is
a subset of the form $p^{-1}(\overline{K})$,
 where $p:\H\ra\R^{2n}$ is the quotient of $\H$
by its center
 and
$\overline{K}\subset\R^{2n}$ is a half-space in the 
usual sense.  (The corresponding result for subsets
 of finite perimeter in $\R^n$ is classical and
 due to DiGiorgi; \cite{degiorgi}.)

Note that for a subset of the form $E=p^{-1}(\overline{K})$,
the associated elementary cut metric
assigns distance $0$ to any pair of points which lie
on the same coset of the center $Z$. In view of (\ref{cutmetrep}),
this strongly suggests that under blow up, a cut metric which is
supported on sets of finite perimeter should
become degenerate in the direction of the center.
Most of our technical 
work consists of making this simple idea rigorous.

\subsection{Metric differentiation and monotonicity}

\mbox{}

\vspace{2mm}

Here we discuss some results related to our main
theorem,  which will appear elsewhere.

There is an alternative approach to the main theorem which
is based on metric differentiation and monotonicity.  We recall
that \cite{kirchheim,pauls} showed that any Lipschitz map
$
f:X\lra Y
$
from a Carnot group $X$ to an arbitrary metric space
$Y$ has a full measure set of points of metric differentiability.
This implies that blow ups of $f$ at almost every point of 
$X$ yield  limit maps 
$
f_\om:X\lra Y_\om,
$
where $Y_\om$ is an ultralimit of rescalings of $Y$.  When $Y=L^1$
then the ultralimit $Y_\om$ is also an $L^1$ space. One can show that 
the cut measure associated
with $f_\om$ is supported by {\em monotone} subsets of $X$; these
are measurable subsets $E\subset X$ such that for almost every 
horizontal geodesic, $\ga\subset X$, the intersection
$E\cap \ga$
is --- modulo a set of zero $1$-dimensional Hausdorff measure ---
either a ray,  the empty set,
or $\ga$ itself.   An analysis of the structure of monotone subsets of the Heisenberg
group eventually leads to another proof of Theorem \ref{thmmainfirst}.  
The details of this will appear in \cite{ckmetricdifferentiation}, 
together with other applications of the same circle of ideas.

Theorem \ref{thmmainfirst} implies that any Lipschitz map from a ball $U\subset\H$
into $L^1$ cannot be bi-Lipschitz.  By a compactness argument, it follows that
if 
$
f:U\ra L^1
$
is an $L$-Lipschitz map, then  the quantity
$$
\eta_L(r)\defeq \inf\left\{\frac{d(f(x),f(y))}{d(x,y)}\;\mid\; x,y\in U,\; d(x,y)\geq r\right\}
$$
can be bounded above by a function $\hat\eta_L:[0,\infty)\ra \R$ which depends
only on $L$, and which satisfies 
$$
\lim_{t\ra 0}\;\hat\eta(t)=0\,.
$$
This leads one to ask for a bound on the asymptotic behavior of $\hat\eta$, which is of
 interest in computer science, in particular, in connection with the 
failure of the Goemans-Linial conjecture.  
Such a bound will be given in a paper with Assaf Naor  \cite{ckn}.

\no
\subsection{Organization  of the paper }

\mbox{}

\vspace{2mm}The remaining sections of the paper are organized as follows.

In Section \ref{prelims}, we collect some 
background material on PI spaces, BV functions
and sets of finite perimeter, which is used
in subsequent sections.

In Section \ref{secl1mappings}, under the additional
assumption that
$X$ carries a measure, $\mu$,
 we give alternative characterizations of 
$L^1$ maps $f:X\to L^1(Y,\nu)$. We also discuss in this
setting, the equivalence
between metrics induced by maps to $L^1$ and cut metrics.

In Section \ref{bv}, assuming in addition that
$X$ is a PI space,
we show the equivalence between metrics induced by
BV maps to $L^1$ and cut measures which are supported
on sets of finite perimeter (FP cut measures). This equivalence
is the basic new conceptual idea in this paper.

In Section \ref{totper}, and for the remainder of the paper,
we consider an FP cut measure $\Si$.
We construct the {\it total perimeter
measure}, $\lambda\in \radon(X)$,
 associated to $\Si$.

In Section \ref{badper}, and in the sections which follow,
 we specialize to the Heisenberg group $\H$.
We specify the bad part of $\lambda$,
 taking into account location and scale.  Here
the ``bad part'' means the part carried by those
cuts which are not close to a half space. Getting
suitable bounds on the bad part of $\lambda$ is
 the key to proving our main differentiation
theorem.

In Section \ref{badpercontrol}, we prove a parameterized
version of the main result of \cite{italians1};
see Theorem \ref{badra0p}. 
This result is of crucial importance, and it
is the only place where we appeal to
\cite{italians1}. From Theorem \ref{badra0p} and
a straightforward argument based on measure 
differentiation, we derive the required 
bounds on the bad part of the perimeter measure.

In Section \ref{gbcuts}, we introduce collections,
$\G,\B$, of good and bad cuts, taking into account
location and scale. Then we translate the estimates of
Section \ref{badpercontrol} into estimates
on $\G$ and $\B$. 

In Section \ref{goodcutmeas}, we construct an FP 
cut measure, $\widehat{\Si}$, associated to $\Si$,
which is supported on cuts which are half spaces.
In constructing $\widehat{\Si}$, we approximate
cuts in $\G$ by cuts which are true half-spaces.

In Section \ref{mainthm},
we prove Theorem \ref{maindiffcut},
our main differentiation theorem. Namely, we show
that at most locations, the normalized
$L^1$-distance between the distance functions induced
by $\Si$ and $\widehat{\Si}$ can be made as small
as we like, provided we go to a sufficiently small
scale. The preceding sections have been organized
in such a way that the proof uses only the estimates
of Sections \ref{gbcuts}, \ref{goodcutmeas}, 
and  the Poincar\'e inequality.

\no
{\bf Acknowledgement.}
We would like to thank Assaf Naor for several crucial contributions
to this work ---  for bringing the Heisenberg embedding 
question to our attention in the first place, for telling us about
Kakutani's theorem \cite{kakutani}, and for drawing our attention to
the computer science literature.

\section{Preliminaries}
\label{prelims}

In this section, we collect some relevant background material
on PI spaces, BV functions and sets of finite perimeter.
\vskip2mm

\no
\subsection{PI spaces.}

\mbox{}

\vspace{2mm}
A {\em PI space} is a  metric measure space, $(X,d^X,\mu)$, for which the metric is complete,
and the doubling condition and Poincar\'e inequality hold.

The doubling
condition on the measure, $\mu$, states that for some $\beta(R)<\infty$,
\begin{equation}
\label{doubling}
\mu(B_{2r}(x))\leq \beta(R)\cdot\mu(B_r(x))\qquad r\leq R\, .
%\mu(B_{2r}(x))\leq 2^{\kappa(R)}\mu(B_r(x))\qquad r\leq R\, .
\end{equation}

A Borel measurable function, $g:X\ra [0,\infty]$, is called an {\it upper gradient}
for $f$ if for every rectifiable curve, 
$c:[0,\ell]\to X$, parameterized by arclength,
$s$,
$$
|f(c(\ell))-f(c(0))|\leq \int_0^\ell g(c(s))\, ds\, .
$$

Put
$$
f_{x,r}=\frac{1}{\mu(B_r(x))}\int_{B_{r}(x)}f\, d\mu\, .
$$

The $(1,1)$-Poincar\'e inequality is the condition
that for some $\tau(R), \lambda<\infty$,
\begin{equation}
\label{pi1}
\int_{B_r(x)}|f-f_{x,\,r}|\,d\mu
\leq r\cdot\tau(R)\int_{B_{\lambda r}(x)}g\,d\mu\, ,
\end{equation}
where $g$ is any upper gradient of $f$.
An equivalent form of the Poincar\'e
inequality is
\begin{equation}
\label{pi2}
\int_{B_r(x)\times B_r(x)}|f(x_1)-f(x_2)|\,d\mu\times\, d\mu
\leq r\cdot\tau'(R)\int_{B_{\lambda r}(x)}g\,d\mu\, .
\end{equation}

 For definiteness and without 
essential loss of generality, in the sequel,
we will assume $\lambda=2$.

 For our present purposes,
it is enough to consider say
$r\leq 1$. Thus, $\kappa(1)=\kappa$, $\tau(1)=\tau$,
$\tau'(1)=\tau'$.

\no
\subsection{The Heisenberg group.}

\mbox{}

\vspace{2mm}The Heisenberg group, $\H$, is a $2$-step nilpotent Lie group
diffeomorphic to $\R^{2n+1}$. When equipped with
the Carnot-Caratheodory metric,
its Hausdorff dimension  $\H$ is $2n+2$.  
We will recall the definition in dimension $3$.
For an extended discussion, see \cite{gromov}. 

The 
 $3$-dimensional
Heisenberg group, $\H\subset\gl(3)$,
consists of matrices,
\begin{equation}
\label{e:1.1}
%  \left\{
\begin{bmatrix}
1 & a&c\\
0&1&b\\
0&0&1
\end{bmatrix}\, ,
\end{equation}
$a,b,c\in\R$.
In particular, $\H$ is diffeomorphic to $\R^3$.

As a vector space, the Lie algebra of $\H$ is 
$\R^{3}=(a,b,c)$, 
realized as the space of matrices,
$$
\begin{bmatrix}
0 & a&c\\
0&0&b\\
0&0&0
\end{bmatrix}\, ,
$$

Let $P=(1,0,0)$, $Q=(0,1,0)$, $Z=(0,0,1)$.

We have the commutation relations for the Lie algebra,
$[P,Q]=Z$, $[P,Z]=[Q,Z]=0$.
In particular, $Z$ is a basis for the center of the Lie algebra
and the center of $\H$ is the $1$-parameter subgroup 
$$
\cent(\H)=\{\exp(tZ)\}_{t\in\R}\, .
$$

To define the Carnot-Caratheodory metric, view 
$P,Q,Z$
as orthonormal left-invariant vector fields on $\H$ and denote by
$\De$, the {\it horizonal distribution},
i.e. the $2$-dimensional distribution on $\H$ spanned by the 
left-invariant vector fields $P,Q$.  
The Carnot-Caratheodory distance, $d^\H(x_1,x_2)$,
between two points $x_1,x_2\in\H$, is defined as the
infimum of the length of paths, 
$\ga:[0,1]\ra\H$, such that $\ga$ 
joins $x_1$ to $x_2$
and the velocity vector of $\ga$ is everywhere tangent to $\De$. 
When equipped with the Carnot-Caratheodory metric and Lebesgue measure,
$\L$ (equivalently, Haar measure), the metric measure
space, $(\H,d^\H,\L)$, is a PI space.
\vskip2mm

\vskip2mm

  \no
 \subsection{Functions of bounded variation on metric measure spaces.}

\mbox{}

\vspace{2mm}

Let $(X,d^X)$, $(W,d^W)$ denote metric spaces.

Given a Lipschitz function, $f\in\Lip(X, W)$, 
the Lipschitz constant, $\LIP\,f$,
is 
$$
\LIP\,f\defeq
\sup_{x,x'}\,\frac{d^W(f(x),f(x'))}{d^X(x,x')}\, .
$$
The 
{\it pointwise Lipschitz constant}, $\Lip\, f$, is 
\begin{equation}
\label{pointwiselip}
\Lip(f(x))
\defeq\liminf_{r\to 0}\,\sup_{d^X(x,x')<r}\frac{d^W(f(x),f(x'))}{r}\, .
\end{equation}
Note that $\Lip\,f$ is an upper gradient for $f$.

Now assume in addition,
that $X$ is a PI space and
let $U\subset X$ denote an open set. Let $V$ denote a Banach space.
\begin{definition}
\label{bvdef}
  The map, $h\in L^1(U,V)$, has {\it bounded variation},
$f\in \bv(U,d^X,\mu,V)$, if there
exists  a sequence of locally Lipschitz functions,
$h_i \stackrel{L_1}{\longrightarrow}h$,  such that 
$$
\liminf_{i\to\infty}\int_U \Lip\, h_i\,d\mu<\infty\, ;
$$
see Definition 2.11 and Remark 2.16 of \cite{cheeger}, and \cite{luigi1,luigi2}.
 \end{definition}

As usual, we just write $f\in\bv(U,V)$ and $f\in\bv(U)$ 
if $V=\R$. 

\begin{remark}
Definition \ref{bvdef} makes sense when the target is an 
arbitrary metric space.
\end{remark}

The {\it variation} of $h\in\bv(U,V)$ is
\begin{equation}
\label{bvug}
\Var(h,U)\defeq\inf_{h_i}\,\liminf_{i\to\infty}\int_U \Lip\, h_i\,d\mu\, .
\end{equation}
When there is no danger of confusion about the domain,
in place of $\Var(h,U)$,
we sometimes write $\Var(h)$.

If $U'\subset U$ is an open subset, there is a natural
variation decreasing restriction map $\bv(U)\to \bv(U')$.
In fact, given  $f\in\bv(U)$,  there is a canonically associated
Radon measure, $\var(h,U)$, on $U$ (sometimes denoted $\var(h)$)
the {\it variation measure} of $h$, whose value on any open
set, $U'\subset U$ is 
\begin{equation}
\label{varmeas}
\Var(h,U')=\var(h,U)(U')\, ;
\end{equation}
see \cite{mir}. 
In particular, for $h\in\bv(U)$,
$$
\Var(h,U)=\mass(\var(h,U))\, ,
$$
where by definition, for $\theta$ a 
measure  on $U$,
\begin{equation}
\label{massdef}
\mass(\theta)=\theta^{pos}(U)-\theta^{neg}(U)\, .
\end{equation}
($\theta^{pos}$, $\theta^{neg}$ denote the 
positive and negative parts of $\theta$, relative 
to the Hahn decomposition.)

The measure,   $\var(h)$,  can also be constructed in a 
manner completely analogous to the construction of
the minimal upper gradient in \cite{cheeger}.
Note that the measure,  $\var(h)$, need {\it not} be absolutely
continuous with respect to $\nu$ e.g. if,
as considered below, $h$ is a characteristic
function $\chi_E$.

It is immediate that the variation is 
lower semicontinuous under $L^1$ convergence. 
The  variation measure  satisfies an analogous
 weak lower semicontinuity  property
under $L^1$ convergence; 
compare Proposition \ref{perinellone}.

By a 
diagonal argument, there 
exists a sequence of locally Lipschitz functions,
$h_i\stackrel{L^1}{\longrightarrow} h$, with
\begin{equation}
\label{relaxed}
\lim_{i\to\infty}\int_U \Lip\, h_i\,d\mu=\var(h)(U)\, .
\end{equation}
Note also that if $U\subset X$, $\mu(U)<\infty$,
and $f:U\to V$ is Lipschitz, then $f\in\bv(U,V)$.

\begin{remark}
 In defining {\it real valued} BV functions on
 $\H$ with its
Carnot-Caratheodory metric, it is equivalent to 
assume $h_i\in C^1$,  and replace $\Lip\, h$, in (\ref{bvug}) by
the norm of the
horizontal derivative --- the restriction of 
the classical differential to the horizontal subspace $\Delta$;
see  \cite{degiorgi,giusti} for the classical
theory of BV functions on $\R^n$ and
\cite{luigi1,luigi2,italians1,italians2}
 for the Heisenberg case.
\end{remark}
\vskip1mm

\begin{remark}
It is clear that if $f\in \bv(X)$ then (\ref{pi1}), (\ref{pi2}),
hold, with the integral of $g$ replaced by 
$\var(f)(B_{\lambda r}(x))$ on the right-hand sides.
\end{remark}
\vskip2mm

\no
\subsection{Sets of finite perimeter.}

\mbox{}

\vspace{2mm} The {\it perimeter measure} of a measurable
set $E\subset U$ is
the variation  measure 
of the  characteristic function $\chi_E$,
\begin{equation}
\label{perdef}
\per(E,U)\defeq \var(\chi_E,U)\, .
\end{equation}
The 
{\it perimeter} of $E$ is the  mass of $\per(E)$,
$$
\Per(E,U)\defeq\mass(\per(E,U))\, .
 $$
The measurable  set $E$ has {\it finite perimeter in $U$} if 
$$
\Per(E,U)=\Var(\chi_E,U)<\infty\, .
$$
As usual, below we tend to supress the dependence on 
$U$.

As above, the perimeter and perimeter measure
are lower semicontinuous (respectively
weakly lower semicontinuous)
under $L^1$
convergence of characteristic functions.

The {\it coarea formula}
for $h\in BV(U)$ functions asserts
\begin{equation}
\label{coarea}
\var(h)(U)=\int_\R\,\Per(\{h\geq t\}\cap U)\,d\L(t)\, ;
\end{equation}
see \cite{ambmirpall}.

\section{$L^1$ maps into $L^1$ spaces}
\label{secl1mappings}

In this section,  $(X,\mu)$, $(Y,\nu)$ will denote $\sigma$-finite measure spaces.

 Here we show that an $L^1$ map
$f:X\to L^1(Y)$ gives rise to an $L^1$ function on the product $X\times Y$,
and an $L^1$ map $g:Y\ra L^1(X)$.   
We also show that the metric, $d_f$, induced by such a map, $f$,
has a cut metric representation, i.e. it is a superposition of
elementary cut metrics,
$$
d_f(x_1,x_2)=\int_{\cut(X)}\; d_E(x_1,x_2)\;d\Si_f(E)\,.
$$
Here $\cut(X)$, $d_E$, and $\Si_f$ are the $L^1$ versions of the 
objects in Subsection \ref{subsecindication}.

We begin with some general remarks and notation.
\vskip2mm

\no
\subsection{$L^1$ maps to Banach spaces.}

\mbox{}

\vspace{2mm}
Denote by $L^1(X,\mu, V)$, the
$L^1$ space of  $(X,\mu)$, with
values in the Banach space $V$. If the
second argument is omitted, we understand
$V=\R$.
Recall that elements of $L^1(X,\mu,V)$ are equivalence classes
of of Borel measurable maps  $f: X\ra V$, for which the norm,
$|f|:X\ra \R$, is an integrable function on $X$.
We will often write $f\in L^1(X,\mu,V)$ when we mean that $f$ is  such an
 equivalence class and refer to the $L^1$ function, $f$,
when we mean that $f$ is a representative of such an equivalence class.

Given $f\in L^1(X,\mu,V)$, there is a well-defined pushforward measure,
$f_*(\mu)$, which is a Borel measure on $V$, with associated
$L^1$-space, $L^1(X,f_*(\mu))$. The induced map (in the opposite direction)
on real valued
functions gives rise to a map, $f^*$ on $L^1$ spaces,
which is an isometric embedding,
$$
f^*:L^1(V,f_*(\mu))\;\lra\;L^1(X,\mu)\, .
$$

We may also use $f$ to pullback the distance function from $V$, 
thereby obtaining
a well-defined equivalence class of measurable  functions, i.e. 
\begin{equation}
\label{eqnmeasdist}
d_f:X\times X\lra \R\, .
\end{equation}
Note that the restriction
\begin{equation}
d_f\restr_{S\times S}:S\times S\lra \R
\end{equation}
is integrable when $\mu(S)<\infty$.

In general,  we use the term {\em $L^1_{\loc}$-distance function} to
refer to equivalence classes of measurable distance functions on $X\times X$
which are integrable on subsets of the form $S\times S$, where $\mu(S)<\infty$.  Note that
it makes sense to integrate a map from a measure space $(Z,\zeta)$
taking values in the space of $L^1_{\loc}$-distance
functions, provided it
becomes an $L^1$ map,
\begin{equation}
\label{eqnl1locintegration}
(Z,\zeta)\lra L^1(S\times S),
\end{equation}
when the distance functions are restricted to $S\times S$, for any finite measure subset
$S\subset X$.
  
\vskip2mm

\no
\subsection{$L^1$ targets; a variant of Fubini's theorem.}

\mbox{}

\vspace{2mm} From now on,  we will usually write $L^1(X)$, $L^1(Y)$ for
$L^1(X,\mu)$, $L^1(Y,\nu)$ respectively and write
$L^1(X\times Y)$ for $L^1(X\times Y,\mu\times\nu)$,
supressing the dependence on the measures.

Given a measurable function $f: X\ra L^1(Y)$ representing an $L^1$ map,
we obtain an element $f(x)\in L^1(Y)$ for each $x\in X$; this is itself
an equivalence class of measurable functions on $Y$.   The main technical point of 
the next result is that one may choose representatives of these equivalence classes
in a measurably varying fashion.

\begin{proposition}
\label{propinterchange}
The spaces
$L^1(X,L^1(Y))$, $L^1(X\times Y)$, and $L^1(Y,L^1(X))$
are canonically isometric. In particular:

\no
{\rm 1)} Given $f\in L^1(X, L^1(Y))$, 
there exists $H\in L^1(X\times Y)$,
such that for a.e. $x\in X$, 
\begin{equation}
\label{eqnF=f}
H(x,y )=f(x)\, 
\quad\mbox{in}\quad L^1(Y)\, .
\end{equation}
Alternatively, by Fubini's theorem, if we view $H$ as an integrable measurable function
on $X\times Y$, then for $\mu\times\nu$ a.e. $(x,y)\in X\times Y$,
\begin{equation}
H(x,y)=f(x)(y).
\end{equation}

\no
{\rm 2)} If $H\in L^1(X\times Y)$,
then for  $\nu$ a.e. $y\in Y$,
the function
$$
g(y)\defeq
\begin{cases}  
H(x,y)\qquad\,\,\,\,\,\, H(x,y)\in L^1(X)\, ,\\

  0\in L^1(X)\qquad{\rm otherwise}\, ,
\end{cases}
$$
defines an element of $L^1(X,L^1(Y))$.
\end{proposition}
\proof
1) By definition, there is a sequence
 of integrable simple maps, 
$$
f_k:X\ra L^1(Y)\, ,
$$
such that
\begin{equation}
\label{eqnconverges}
\lim_{k\to\infty}\|f-f_k\|_{L^1}= 0\, .
\end{equation}
Without loss of generality, we may assume that  in addition,
the sequence has bounded variation in $L^1(X,L^1(Y))$, i.e.
\begin{equation}
\label{eqnsummable}
\sum_k\|f_{k+1}-f_k\|_{L^1}<\infty\, .
\end{equation}

 For each $k$, the map $f_k$ takes finitely many values;
for each of these 
we pick a measurable representing function, and thereby get
a function 
$H_k:X\times Y\ra \R$, which is clearly measurable.  By 
(\ref{eqnsummable}), for $\mu$ a.e. $x\in X$,
the sequence of integrable functions, $H_k(x,\cdot)$, has
bounded variation in $L^1(Y)$:
$$
\sum_k 
\|H_{k+1}(x,\, \cdot\, )-H_k(x,\, \cdot\, )\|_{L^1}<\infty\, .
$$
Therefore, the sequence, $H_k$, converges pointwise $\mu\times\nu$ almost everywhere.
Thus, we get a measurable function,
$$
H\defeq \liminf_{k\ra\infty} H_k\, ,
$$
which is integrable by Fubini's theorem, and as a consequence of
(\ref{eqnconverges}),
 satisfies (\ref{eqnF=f}).

\medskip
2) This follows by approximating 
the positive (respectively negative) part of 
$H$ 
by a monotone nondecreasing (respectively decreasing) 
sequence of
 functions $H_k$,
where each $H_k$ is a finite linear combination 
of characteristic functions
of rectangles in $X\times Y$.

It is clear that the constructions in 1) and 2) above define isometries
which, by Fubini's theorem, are inverses of one another.
\qed
\vskip2mm

\no
\subsection{Borel measures on $L^1$ and tautological maps.}

\mbox{}

\vspace{2mm}By Proposition \ref{propinterchange}, an $L^1$ map,
$f:X\ra L^1(Y)$,  
induces an $L^1$ map $g:Y\ra L^1(X)$.  

\begin{definition}
\label{deft_f}
The Borel measure, $\T_f$, on  $L^1(X)$ is the measure
$$
\T_f\,:=\, g_*(\nu)\, .
$$ 
\end{definition}

By Fubini's theorem,
we have
\begin{equation}
\label{eqncutmeasure1}
\int_{L^1(X)}\;\|u\|_{L^1}\,d\T_f(u)=\|f\|_{L^1}=\|g\|_{L^1}<\infty\, .
\end{equation}

More generally, let $\T$ denote an
arbitrary Borel measure on $L^1(X)$ satisfying the integrability
condition
\begin{equation}
\label{eqnfinite}
\int_{L^1(X)}\|u\|_{L^1}\, d\T(u)<\infty\, .
\end{equation}

The identity map 
$$
(L^1(X),\T)\lra L^1(X,\mu),
$$
where the domain is viewed as a measure space, 
and the target is viewed as an $L^1$ space, satisfies the 
hypotheses of Proposition \ref{propinterchange}, where
$(L^1(X),\T)$ plays the role of $(X,\mu)$, and $L^1(X,\mu)$
plays the role of $L^1(Y,\nu)$.  This yields:

\begin{corollary}
\label{cortautological}
There is an $L^1$  function 
\begin{equation}
\La\in L^1(L^1(X)\times X,\T\times\mu)
\end{equation}
and an $L^1$ map
\begin{equation}
\taut_{\T}:X\ra L^1(L^1(X),\T)
\end{equation}
such that: 

\no
${\rm 1})$  For any representative of $\La$, we have 
\begin{equation}
\La(u,x)=u(x)\quad\mbox{for $\T\times\mu$ a.e. $(u,x)\in L^1(X)\times X$}.
\end{equation}

\no
$\rm{2})$ For any representative of $\La$, we obtain a representative
of $\taut_\T$ by the formula
\begin{equation}
\taut_\T(x)= \La(\,\cdot\,,x).
\end{equation}

Note that 
$$
\|\taut_\T\|_{L^1}=\int_{L^1(X)}\,\|u\|_{L^1}d\T(u).
$$
\end{corollary}
\proof
Observe that $\T$ is a $\si$-finite measure, since the function
$$
\|\cdot\|:L^1(X)\ra\R
$$
is integrable with respect to $\T$.  Also, the identity map
$$
L^1(X,\mu)\lra L^1(X,\mu)
$$
is Borel measurable, and by (\ref{eqnfinite}), determines
an $L^1$ map
$$
(L^1(X),\T)\lra L^1(X).
$$
We now apply Proposition \ref{propinterchange} with
$$
f=\id_{L^1(X)},
$$
and set $\La\defeq H$.  The lemma  follows.

\qed

\begin{lemma}
\label{lemtautt_f}
Let $f:X\ra L^1(Y)$ be an $L^1$ map, and 
$$
\taut_{\T_f}:X\lra L^1(L^1(X),\T_f)
$$ 
be the map of Corollary \ref{cortautological}.
Then 

\no
$\rm{1})$
$f=g^*\circ \taut_{\T_f}$, where
$$
X\stackrel{\taut_{\T_f}}{\lra}\,L^1(L^1(X),\T_f)\stackrel{g^*}{\lra} L^1(Y,\nu).
$$

\no
$\rm{2})$ The distance functions induced by $f$ and $\taut_{\T_f}$
coincide.  (Recall that as in (\ref{eqnmeasdist}) these are equivalence
classes of measurable functions on $X\times X$.)
\end{lemma}
\proof
Let $\La$ be as in Corollary \ref{cortautological},
so the map 
$$
x\mapsto \La(\,\cdot\, ,x)
$$
is a representative of $\taut_{\T_f}$,
and
\begin{equation}
\label{eqnLaprop}
\La(u,x)=u(x)\quad\mbox{for}\quad \T_f\times\mu\;\;\mbox{a.e.}\quad (u,x)\in L^1\times X.
\end{equation}
Let $H:X\times Y\ra\R$
be as in Proposition \ref{propinterchange}, so 

\begin{equation}
\label{eqnHprop}
H(x,y)=g(y)(x)=f(x)(y)
\end{equation}
for $\mu\times\nu$ a.e. $(x,y)\in X\times Y$.
Since $\T_f=g_*(\nu)$, (\ref{eqnLaprop}) and (\ref{eqnHprop}) imply that

\begin{equation}
\label{eqnLaH}
\La(H(\, \cdot\, ,y),x)=H(x,y)
\end{equation}
for $\mu$ a.e. $x$ and $\nu$ a.e. $y$.
Therefore for $\mu$ a.e. $x\in X$ and $\nu$ a.e. $y\in Y$, 
\begin{equation}
\begin{aligned}
\left(g^*\circ\taut_{\T_f}(x)\right)(y)&=g^*(\La(\,\cdot\,,x))\;(y)\\
&=\La(H(\,\cdot\,,y),x)\\
&=H(x,y)\quad\mbox{by (\ref{eqnLaH})}\\
&=f(x)(y),
\end{aligned}
\end{equation}
which implies 1).

Assertion 2) follows immediately from 1)  because $g^*$ is an isometric embedding.

\qed
\vskip2mm

\no
\subsection{Cut measures}

\mbox{}

\vspace{2mm}\begin{definition}
A {\em  cut } in $X$ is an equivalence class of finite measure
subsets of $X$.  
\end{definition}

We denote the set of  cuts in $X$ by $\cut(X)$,
and  identify it with the set of elements of $L^1(X)$
which can be represented by
characteristic functions.  This 
is a closed
subset of $L^1(X)$. In particular, $\cut(X)$ inherits a metric from
$L^1(X)$.

\begin{definition}
\label{defcutmeasuremetric}
A Borel measure $\Si$ on $\cut(X)$ is a {\em  cut measure} if
$$
\int_{\cut(X)}\;\|\cdot\|_{L^1}\;d\Si\;<\;\infty.
$$
\end{definition}

Since $\cut(X)$ is a closed subset of $L^1(X)$, we may view $\Si$ as 
a measure satisfying (\ref{eqnfinite}).  Therefore by Corollary \ref{cortautological} 
we obtain a tautological map
\begin{equation}
\label{eqntautsi}
\taut_{\Si}:X\ra L^1(\cut(X),\Si),
\end{equation}
where we have used the fact that 
$L^1(L^1(X),\Si)$
is isometric to $L^1(\cut(X),\Si)$.

  Next, using slices, we show how a measure $\T$ 
on $L^1(X)$ satisfying (\ref{eqnfinite})  
 gives rise
to an associated cut measure $\Si_{\T}$. 
 
\bigskip
\begin{lemma}
Let $\super$ be the map 
\begin{equation}
\super:L^1(X)\times \R\lra \cut(X)
\end{equation}
be given by 
\begin{equation}
\super(u,t)\defeq 
\begin{cases}
\quad \{u\geq t\} & \text{when $t>0$}\,,\\
\quad \emptyset& \text{when $t=0$}\,,\\
\quad \{u\leq t\} & \text{when $t<0$\,.}
\end{cases}
\end{equation}
Then

\no
$\rm{1})$ $\super$ is well-defined.

\no
$\rm{2})$ $\super$ has a set-theoretic semicontinuity property: if $(u_k,t_k)\in L^1(X)\times \R$
is a sequence converging to $(u,t)$, then
\begin{equation}
\mu(\super(u_k\,,\, t_k)\setminus \super(u\,, t))\ra 0\;\mbox{as}\;k\ra\infty.
\end{equation}

\no
$\rm{3})$ $\super$ is Borel measurable.

\end{lemma}
\proof
It suffices to consider the case when $t\geq 0$, and the functions 
are nonnegative.

The map $\super$ is well-defined, because if
  two measurable functions $u,v$ represent the same element of 
$L^1(X)$, then for every $t\in \R$, the symmetric difference
\begin{equation}
\{u\geq t\}\;\De\;\{v\geq t\}
\end{equation}
has measure zero, and hence the two sets determine the same element of 
$\cut(X)$.

We now prove 2).  Pick $\de>0$.  Then 
\begin{equation}
\|u_k-u\|_{L^1}=\int_X|u_k-u|d\mu\geq (\de+t_k-t)\;\mu(\{u_k\geq t_k\}\setminus
\{u\geq t-\de\}),
\end{equation}
which forces
\begin{equation}
\mu(\{u_k\geq t_k\}\setminus\{u\geq t-\de\})\ra 0\quad\mbox{as}\quad k\ra\infty.
\end{equation}
Since
\begin{equation}
\mu(\{u\geq t-\de\}\setminus\{u\geq t\})\ra 0\quad\mbox{as}\quad \de\ra 0,
\end{equation}
 this implies 2).

Borel measurability of $\super$ follows from the fact that the collection of
open sets
$$
U(E,\eps)\defeq \{E'\in\cut(X)\mid \mu(E'\setminus E)<\eps\}
$$
generates the full Borel $\si$-algebra, and by assertion 2), 
$\super^{-1}(U(E,r))$ is open in $L^1(X)\times\R$ for all $E\in \cut(X),\;r>0$.

\qed

\bigskip
Given a Borel measure $\T$ on $L^1(X)$ satisfying (\ref{eqnfinite}),
we obtain a cut measure 
\begin{equation}
\label{eqnsi_t}
\Si_{\T}=\super_*(\T\times\L)\,.
\end{equation}

\begin{definition}
\label{deffcutmeasure}
The  {\em  cut measure associated with an $L^1$ map $f:X\ra L^1(Y)$} 
is the Borel measure $\Si_{\T_f}$,
where $\T_f$ is as in (\ref{deft_f}).
\end{definition}

\subsection{The cut metric representation}
\label{seccutmetricrepresentation}

\begin{definition}
\label{defcutmetric}
The {\em elementary  cut metric  $d_E$} associated with a cut 
$E\in \cut(X)$ is the $L^1_{\loc}$-distance function given by 
$$
d_E(x_1,x_2)=|\chi_E(x_1)-\chi_E(x_2)|.
$$

The {\em  cut metric} $d_\Si$ associated with a cut measure $\Si$ 
is the corresponding superposition of elementary cut metrics:
\begin{equation}
\label{eqncutmetric}
d_{\Si}\;=\;\int_{\cut(X)}\;d_E\;d\Si(E).
\end{equation}
Here we view the integration on the right hand side as taking place in the space
of $L^1_{\loc}$-distance functions,
as in (\ref{eqnl1locintegration}).
\end{definition}

\begin{proposition}
\label{propdsi=rhosi}
Let $\Si$ be a cut measure, and $\T$ be a measure satisfying (\ref{eqnfinite}).
 Then: 

\no $\rm{1})$
 The  distance function induced
by the tautological map
$$
\taut_{\Si}:X\lra L^1(L^1(X),\Si)
$$
coincides with the cut metric $d_{\Si}$.

\no
$\rm{2})$ If $\Si=\Si_\T$, then
$$
d_{\taut_{\T}}=d_{\taut_{\Si}}=d_{\Si}.
$$
\end{proposition}
\proof
Let 
$$
\La_{\T}:L^1(X)\times X\lra \R,\quad \La_{\Si}:\cut(X)\times X\lra \R
$$
be  measurable functions as in Corollary \ref{cortautological}.  Then
for $\mu\times \mu$ a.e. $(x_1,x_2)\in X\times X$,

\begin{equation}
\begin{aligned}
d_{\taut_{\Si}}(x_1,x_2)=&\int_{\cut(X)}\,\left|
\La_{\Si}(E,x_1)-\La_{\Si}(E,x_2)\right|d\Si(E)\\
=&\int_{\cut(X)}|\chi_E(x_1)-\chi_E(x_2)|\;d\Si(E)\\
=&\int_{\cut(X)}d_E\;d\Si(E)\\
=&\;d_{\Si}\,.
\end{aligned}
\end{equation}

If $\Si=\Si_\T$, then by the definition of $\Si_\T$ as the pushforward
of $\T\times\R$ under $\super$, we may continue the calculation:

\begin{equation}
\begin{aligned}
d_{\taut_{\Si}}=&\int_{\cut(X)}|\chi_E(x_1)-\chi_E(x_2)|\;d\Si(E)\\
=&\int_{L^1(X)\times\R}|\chi_{\super(u,t)}(x_1)-\chi_{\super(u,t)}(x_2)| d(\T\times\L)(u,t)\\
=&\int_{L^1(X)}\int_\R |\chi_{\super(u,t)}(x_1)-\chi_{\super(u,t)}(x_2)|d\L(t)\,d\T(u)\\
=&\int_{L^1(X)}|u(x_1)-u(x_2)|d\T(u)\\
=&d_{\taut_{\Si_\T}}.
\end{aligned}
\end{equation}

\qed

\begin{proposition}
The distance function induced by an $L^1$ map
$$
f:X\lra L^1(Y)
$$ 
is the same as that
induced by the tautological map
$$
\taut_{\Si_f}:X\lra L^1(\cut(X),\Si_f).
$$
\end{proposition}
\proof
By Lemma \ref{lemtautt_f} and Proposition \ref{propdsi=rhosi}, we have
\begin{equation}
d_f=d_{\taut_{\T_f}}=d_{\taut_{\Si_f}}\,.
\end{equation}

\qed

\section{BV maps to $L^1$ and FP cut measures}
\label{bv}

  We retain  our notation from  the previous section.  Thus $(X,\mu)$ and $(Y,\nu)$ 
will be  $\si$-finite measure spaces. However,
 we assume in addition that
$X$ carries a metric,
$d^X$,   such that $(X,d^X,\mu)$ is a PI space, i.e.
 $\mu$ satisfies
a doubling condition and 
a $(1,1)$-Poincar\'e inequality;  we let $\kappa$ and $\tau$
be as in Section \ref{prelims}.

The key new observation of this paper can be summarized as follows.

Suppose  $f:X\ra L^1(Y)$ is a map of bounded variation; for instance $f$ could
be any   Lipschitz map, provided $\mu(X)<\infty$.  Let $g:Y\ra L^1(X)$ be the $L^1$ map
provided by Proposition \ref{propinterchange}.  Now,
 the roles of $X,f$ and $Y,g$ are no longer symmetrical.

Although the regularity of the map $g$ is worse than that of $f$ --- it is typically
only measurable whereas $f$ is BV --- the typical function, $g(y)\in L^1(X)$,  
has better regularity than the typical function $f(x)\in L^1(Y)$:
$g(y)$ has {\em bounded variation}, $\Var(g(y))<\infty$,
and the the integral over $Y$  of the function,
$\Var(g(y))$, is {\it finite}. 
In fact, these conditions provide 
a characterization of BV maps to $L^1$; see 
Theorem \ref{lembvslices}.

We also give a second and, in a sense, more directly relevant
characterization of BV maps to $L^1$, in terms of 
what we call 
``FP cut measures'' (where FP stands for finite perimeter).
We show that $f\in\bv(U,L^1(Y))$ if and only if the cut measure,
$\Si_f$, is an FP cut measure. 
Essentially, this follows from the previous 
characterization via the coarea formula.

\begin{remark}
\label{kembed}
By way of contrast with the case of $L^1$ targets, 
note that for the Kuratowski embedding
of $(X,d^X)$ 
into $L^\infty(X,d^X)$, we have $X=Y$, $f=g$, and
nothing is gained. On the other hand, our present point of view
may be useful when studying other function space targets.
\end{remark}
\vskip2mm

\no
\subsection{Characterizing BV maps to $L^1$  by variation.}

\mbox{}

\vspace{2mm}Let $U\subset X$ denote an open subset.
Let $f\in L^1(U, L^1(Y))$ and let $H$, $g$ denote
the  maps in  Proposition \ref{propinterchange}.

 Note that since $\Var(\,\cdot\, )$ 
is a lower semicontinuous function on $L^1(U)$, the integral
$$\int_Y\Var(g(y),U)\, d\nu$$
is a well-defined  extended real number. 

\begin{definition}
\label{finitetotvar}
The map, $f\in L^1(U,L^1(Y))$ has {\it finite total variation},
if $g(y)\in\bv(U)$, for
$\nu$ a.e. $y\in Y$ and
\begin{equation}
\label{aebv}
\int_Y\Var(g(y),U)\, d\nu<\infty\, . 
\end{equation}
The quantity in (\ref{aebv}) is the
{\it total variation} of $f$.
\end{definition}
\vskip2mm

The following theorem shows that the total variation, which is defined
only for $L^1$ targets, is comparable to the variation defined in (\ref{bvug}),
which is defined for arbitrary Banach space targets.

\begin{theorem}
\label{lembvslices}
 $f\in \bv(U,L^1(Y))$ if and only if 
$f$ has finite total variation.
Moreover, there is a constant, $c=c(\kappa,\tau)>0$,
such that
\begin{equation}
\label{uppergradbound}
c^{-1}\cdot\Var(f,U)\leq\int_Y\Var(g(y),U)\,d\nu\leq c\cdot\Var(f,U)\, .
\end{equation}
\end{theorem}
 
\proof
Assume $f\in\bv(U,L^1(Y))$.  
Since $\Var(g(y),U)<\infty$ implies that
$\var(g(y),U)$ is Borel regular, 
by the monotone convergence theorem,
it suffices to consider
an open set
 $U'\subset U$ with compact closure in $U$, and
to establish the inequalities in (\ref{uppergradbound})
with $\Var(g(y),U)$ calculated on $U'$ rather than on $U$.

By (\ref{relaxed}),  there exists a sequence,
$f_i\in \Lip(U,L^1(Y))$, with
$f_i \stackrel{L_1}{\longrightarrow}f$, such that
$$
\lim_{i\to\infty}\int_U \Lip\,  f_i\,d\mu=\Var(f,U)\, .
$$

 Fix an open set $U'\subset U$ with compact closure in $U$.  
We will construct a sequence,
$f_{i,j}\in \Lip(U',L^1(Y))$,
with
$f_{i,j}\stackrel{C^0}{\longrightarrow}f_i$ on $U'$,
such that for $c=c(\kappa,\tau)$, 
\begin{equation}
\label{varij}
\int_{U'} \Lip\, f_{i,j}\,d\mu
\leq c\cdot\int_{U'} \Lip\, f_i\,d\mu\, ,
\end{equation}
\begin{equation}
\label{varijs}
\int_Y \left(\int_{U'} \Lip\, g_{i,j}(x,y)\,d\mu\right)d\nu\,
\leq c\cdot\int_{U'} \Lip\, f_{i,j}\,d\mu\, .
\end{equation}
Then the claim follows by a diagonal argument,
together with Fatou's lemma.

Let  $\{x_{i,j,k}\}$ denote a maximal $j^{-1}$ separated
set in $U$. By a standard lemma, the multiplicity of
the covering, $\{B_{2j^{-1}}(x_{i,j,k})\}$, is bounded by
$N=N(\kappa)$. Also, by using distance functions
from the points, $x_{i,j,k}$, 
we can construct in standard fashion, a 
partition of unity, $\{\phi_{i,j,k}\}$, 
subordinate to $\{B_{2j^{-1}}(x_{i,j,k})\}$, with

\begin{equation}
\label{LIPphi}
 \LIP(\phi_{i,j,k})\leq c(\kappa)\cdot j\, .
\end{equation}

Define the {\it regularization}, $f_{i,j,k}$, of $f_{i,j}$, by
\begin{equation}
\label{reg}
f_{i,j}=\sum_k \,\overline{f}_{i,j,k}\cdot \phi_{i,j,k}\, ,
\end{equation}
where
$$
\overline{f}_{i,j,k}=\frac{1}{\mu(B_{2j^{-1}}(x_{i,j,k}))}
\int_{B_{2j^{-1}}(x_{i,j,k})}\, d\mu\,  .
$$

 Since, $f_i$ is Lipschitz, it follows that
$f_{i,j}\stackrel{C^0}{\longrightarrow}f_i$.

 From now on, we only consider $j$ so large that 
if $\supp(\phi_{i,j,k})\cap U'\ne \emptyset$,
then $B_{8j^{-1}}(x_{i,j,k})\subset U$.

Let $\ell$
denote a linear functional of norm $1$ on $L^1(Y)$.
Then
$$
\ell(\overline{f}_{i,j,k})
=\int_{B_{2j^{-1}}(x_{i,j,k})}\ell(f_i)\,d\mu\, .
$$
By applying the 
Poincar\'e inequality on $B_{8j^{-1}}(x_{i,j,k})\subset U$ to
the Lipschitz function $\ell\circ f_i$ for all
such $\ell$, and using
 the Hahn-Banach theorem, we conclude that
for all $x_{i,j,k}$, $x_{i,j,k'}$, with
$d^X(x_{i,j,k},x_{i,j,k'})\leq 4j^{-1}$,
\begin{equation}
\label{ugbound}
\|\overline{f}_{i,j,k}-\overline{f}_{i,j,k'}\|_{L^1}
\leq c\cdot j^{-1} \frac{1}{\mu(B_{8j^{-1}}(x_{i,j,k}))}
\cdot\int_{B_{8j^{-1}}(x_{i,j,k})}\Lip\, f_i\,d\mu\, ,
\end{equation}
where $c=c(\kappa,\tau)$.

For any fixed index, $k^*$, we can write
\begin{equation}
\label{diff}
f_{i,j}=\overline{f}_{i,j,k^*}+
        \sum_k\, (\overline{f}_{i,j,k}-\overline{f}_{i,j,k^*})
                  \cdot\phi_{i,j,k}\,.
\end{equation}
Since, $\Lip\, f_i\leq\LIP\, f_i<\infty$, from 
(\ref{LIPphi}), (\ref{ugbound}), and the Lebesgue differentiation
theorem applied to $\Lip\, f_i$, 
 we easily get (\ref{varij}).

Since $\overline{f}_{i,j,k}\in L^1(Y)$,
relation (\ref{varijs}) follows from (\ref{diff}), 
and a straightforward argument based on Fubini's theorem.

 Now assume, conversely, that 
$g(y)\in\bv(U)$, for
$\nu$ a.e. $y\in Y$ and
that (\ref{aebv}) holds.

By an exhaustion argument it is easily
checked that it suffices to assume that $\nu(Y)<\infty$.
Similarly, by a truncation argument, one can
assume that $H(x,y)$ is bounded.

Define the regularization, $f_j$, of $f$
as in (\ref{reg}). Then $f_j$ and $H_j$ are Lipschitz
and for $\nu$ a.e. $y\in Y$,
the function, $H_j$, is equal to the corresponding
regularization $g_j(y)$ of $g(y)$. 
Moreover, $f_j\stackrel{L^1}{\longrightarrow} f$.

By arguing as above (compare the
verification of (\ref{varij}), (\ref{varijs}))
 and using Fubini's theorem, we have
$$
\begin{aligned}
\infty &>c\cdot\int_Y\var(g(y),U)\, d\nu\\
       &\geq \int_Y\left(\int_{U'} \Lip\,H_j(x,y)\, d\mu\right)\, d\nu\\
       &=\int_{U'}\left(\int_Y \Lip\,H_j(x,y)\, d\nu\right)\, d\mu\\
       &\geq\int_{U'} \Lip\, f_k\, .
\end{aligned}
$$
This suffices to complete the proof.
\qed

\begin{remark}
In actuality, a metric measure space with the doubling property
satisfies a Poincar\'e inequality for real valued functions
if and only if it satisfies a Poincar\'e inequality for
functions with values in an arbitrary Banach space;
see \cite{HeKoShTy}. 
In justifying (\ref{ugbound}) above, rather than using this result,
we appealed directly to the Hahn-Banach theorem.
\end{remark}

 \no
 \subsection{BV maps to $L^1$ and FP cut measures.}

\mbox{}

\vspace{2mm}Let 
$$
\Per_U:\cut(U)\to[0,\infty]\, ,
$$ 
be given by
$$
E\mapsto \Per(E,U)\, .
$$
\begin{definition}
\label{deffiniteperimetercutmeasurefp}
A cut measure $\Si$ is an 
 {\it FP cut measure}
if $\Per_U\in L^1(\cut(U),\Si)$: 
\begin{equation}
\label{FPcutmeas}
\int_{\cut(U)}\Per(E,U)\, d\Si<\infty.
\end{equation}
The quantity in (\ref{FPcutmeas}) is the
{\it total perimeter} of $\Si$.
\end{definition}

\begin{definition}
\label{deffiniteperimetercutmeasure}
$E\in \cut(U)$ is an {\it FP cut} if
$\Per(E,U)<\infty$. 
\end{definition}

Let $\fp(X)\subset \cut(U)$ denote the collection
of FP cuts.
Since $\Per_U\equiv\infty$
on $\cut(U)\setminus\fp(U)$, it follows that
$$
\Si\left(\cut(U)\setminus\fp(U)\right)=0\, .
 $$ 

Let $\taut_{\Si}:U\to L^1(\cut(U),\Si))$ be as in (\ref{eqntautsi}).

\bigskip
\begin{proposition}
A cut measure, $\Si$, is an FP cut measure
if and only if 
 $\taut_{\Si}\in\bv(U,L^1(\cut(U),\Si))$.
\end{proposition}
\proof
The map, $g:\cut(U)\ra L^1(U)$,
 associated with $\taut_{\Si}:U\ra L^1(\cut(X),\Si)$
is given by $g(E)=\chi_E$.  Therefore, by Theorem \ref{lembvslices},
the map $\taut_{\Si}$ is $\bv$ if and only if
$$
\infty>\int_{\cut(U)}\Var(\chi_E)\, d\Si
=\int_{\cut(U)}\Per(E,U)\, d\Si\, ,
$$
if and only if $\Si$ is an FP cut measure.
\qed
\vskip2mm

The next proposition asserts the equality of the
total variation of $f$ and total perimeter of
$\Si_f$.
\begin{proposition}
\label{bvcutchar}
If $f\in L^1(U,L^1(Y))$,  then 
$\Si_f$ is an FP cut measure
if and only if $f\in\bv(U,L^1(Y))$. Moreover,
\begin{equation}
\label{totperequalstotvar}
\int_{\cut(U)}\Per(E,U)\, d\Si_f
=\int_Y\Var(g(y,U))\, d\nu\, .
\end{equation}
\end{proposition}
\proof

We define 
$$
S:\R\times Y\lra \cut(X)
$$
by 
$$
S(t,y)\defeq \super(g(y),t),
$$
 where 
$$
g:Y\lra L^1(U)
$$
is the map of Proposition \ref{propinterchange}.
By  the definition
of the cut measure $\Si_f$, Fubini's theorem, and (\ref{coarea}), we have

\begin{equation}
\begin{aligned}
\int_{\cut(U)}\Per(E,U)\, d\Si_f=& \int_{Y\times\R} \Per(S(t,y),U)
\, d(\L\times\nu)\\
=&\int_Y\,\int_\R \Per(S(t,y),U)\, d\nu\, d\L\\
=&\int_Y\Var(g(y),U)\, d\nu \,.
\end{aligned}
\end{equation}
 Therefore by Theorem \ref{lembvslices}, the map $f$
is $\bv$ if and only if $\Si$ is an FP cut measure.
\qed

\section{The total perimeter measure}
\label{totper}

We retain the notation of the preceding section.

In this section we will associate to each FP
cut measure, $\Si$, a 
Radon measure, $\lambda_\Si\in\radon(U)$, called 
the  total perimeter measure of $\Si$, whose   
mass  is the total perimeter of $\Si$.
The measure $\lambda_\Si$ is obtained by 
integrating
the measure valued function,
$$
\per(E,U):\cut(U)\to\radon(U),
$$ with respect to
$\Si$.

In the main result on the Heisenberg group,
an essential point is to suitably control the 
``bad part'' of $\lambda_\Si$; see 
Sections \ref{badper}--\ref{mainthm}.

\subsection{Integrating measure valued functions.}

\mbox{}

\vspace{2mm}Let $(Z,\zeta)$ denote a measure space.

Let  $L$ denote a locally
compact Hausdorff space and $C_c(L)$  the space of continuous
functions of compact support, equipped with the sup norm.

A map,
$$
\Psi:(Z,\zeta)\ra \radon(L)\, ,
$$
is {\it weakly measurable}  if for every $\phi\in C_c(L)$,
\begin{equation}
\label{Psi}
z\mapsto \int_L \phi\; d\Psi
\end{equation}
is a measurable function on $Z$.

The map $\Psi$ is
  {\it weakly $L^1$} if it is weakly measurable and there
exists $C<\infty$, such that for all $\phi\in C_c(L)$,
\begin{equation}
\label{Psiellone}
\int_Z\int_L \phi \,\, d\Psi \,d\zeta
\leq C\cdot \|\phi\|_{L^\infty}\, .
\end{equation}

According to the next proposition, 
a weakly $L^1$ map into $\radon(L)$ can be integrated to obtain a Radon
measure.
\begin{proposition}
\label{lemtotalmeasure}
Let $\Psi: (Z,\zeta)\ra \radon(L)$ denote a weakly $L^1$ map. Then
there is a  measure, $\eta\in \radon(L)$, such that for every
Borel set $A\subset L$, 
$$
\eta(A)=\int_Z\Psi(z)(A)\, d\zeta\, .
$$

If the measure $\Psi(z)$ is nonnegative
for $\zeta$ a.e. $z\in Z$,
then
\begin{equation}
\label{eqntotalmass}
 \mass(\eta)=\int_Z\mass(\Psi(z))\, d\zeta.
\end{equation}
\end{proposition}
\proof
Since $\Psi$ is weakly $L^1$, it follows that 
the formula
\begin{equation}
\phi \mapsto\int_Z\left( \;\int_L\phi(x)\;  d(\Psi(z))(x)\right)\; d\zeta(z)\, .
\end{equation}
defines a bounded linear functional on $C_c(L)$.
Thus, the proposition  follows from the Riesz representation theorem.
\qed  
\vskip2mm

\no
\subsection{Constructing the total perimeter measure $\lambda_\Si$.}

\mbox{}

\vspace{2mm}
\begin{proposition}
\label{perinellone}
 Given an FP cut measure $\Si$, 
the map given by
$$
E\mapsto\per(E,U)
$$
defines a weakly $L^1$ map,
$$
(\cut(U),\Si)\to \radon(U)\, .
$$
\end{proposition}
\proof
By essentially the same observation as that which shows
that $\per(E,U)$ is lower semicontinous under
$L^1$ convergence of characteristic functions,
it follows that the map in (\ref{Psi})
is the difference of two lower semicontinuous
functions (corresponding to the nonnegative
and nonpositive parts of the function $\phi$).
It is then clear that (\ref{Psiellone}) holds.
\qed

\bigskip
\begin{definition}
\label{totalper}
The {\it total perimeter measure}
$\lambda_\Si\in \radon(U)$
of the FP cut measure
$\Si$
is the measure obtained by integrating
the weakly $L^1$ map 
$E\to \per(E,U)$.
\end{definition}

\bigskip

\begin{remark}
Note that by (\ref{eqntotalmass}), 
\begin{equation}
\mass(\lambda_\Si)=\int_{\cut(U)}\Per(E,U)\,d\Si\, ;
\end{equation}
see Definition \ref{deffiniteperimetercutmeasurefp}.
In case $\Si=\Si_f$ for some some
$f\in\bv(U,L^1(Y))$, the total perimeter of 
$\Si_f$ is equal to the total variation of $f$;
see Definition \ref{finitetotvar} and (\ref{totperequalstotvar})
of Proposition \ref{bvcutchar}.
\end{remark}

\no
\subsection{Lipschitz maps to $L^1$.}

\mbox{}

\vspace{2mm}\begin{proposition}
\label{abcon}
There is a constant, $C<\infty$, depending only on the
constants $\be$, $\la$  and $\tau$, with the following property.
A BV map $f:U\ra L^1(Y)$ admits 
an $L$-Lipschitz representative  
if and only if for every ball, $B_r(x)$,
\begin{equation}
\label{lipineq}
\frac{\lambda_{\Si_f}(B_r(x))}{\mu(B_r(x))}\leq C \cdot L\, .
\end{equation}
\end{proposition}
\proof
Since, Proposition \ref{abcon} is not required in the sequel,
we will be very brief.  

The necessity of (\ref{lipineq}) follows from the argument used
in proving Theorem \ref{lembvslices}.

The sufficiency follows
from an application of the ``telescope estimate'' as in the
proof of the standard estimate, (4.19), of \cite{cheeger}.
Here is a sketch of a variant of that argument.
One considers a pair of points $x,x'\in U$,
and for small $r>0$, a suitably chosen sequence of points
$$
x=x_1,\ldots,x_k=x'
$$
where $d(x_i,x_{i+1})<\frac{r}{2}$, and 
$$
k\leq \const \frac{d(x,x')}{r}.
$$
The Poincar\'e inequality and (\ref{lipineq}) imply that there is a constant $C=C(\be,\la,\tau)$, 
such that for all $1\leq i<k$, the average
of $f$ over $B_r(x_i)$ differs by at most $C L r$
from its average over $B_r(x_{i+1})$.  So 
$$
\left|\;\frac{1}{\mu(B_r(x))}\int_{B_r(x)} f\,d\mu-
\frac{1}{\mu(B_r(x))}\int_{B_r(x)} f\,d\mu\;\right|
$$
$$
\leq CLkr<C' L d(x,x'),
$$
where $C'=C'(\be,\la,\tau)$.
Since this estimate is independent of $r$, 
it follows that $f$ has a $C'L$-Lipschitz representative. 
\qed

\section{The total bad perimeter measure}
\label{badper}
 
We retain the notation of the preceding sections, except
that we will just write $\per(E)$ in place of 
$\per(E,U)$, supressing the dependence on $U$.

In the remaining sections, we are concerned with 
properties of
sets of finite perimeter which are not valid in
general PI spaces.  For this reason, from now on, $X$ will
be either $\R^n$ or the Heisenberg group $\H$
with its Carnot-Caratheodory metric, 
$\mu$ will denote Lebesgue measure (or equivalently
Hausdorff measure) and $U$ will 
denote a ball in $X$.
In actuality, the discussion has a direct extension 
to the case in which $X$ is replaced
by any 2-step nilpotent Lie group.

We call  a subset
 $E\subset X$ a {\em half-space} if either $X=\R^n$
and $E$ is a half-space in the usual sense, or $X=\H$
and $E$ is the inverse image of a 
Euclidean
half-space under the quotient homomorphism, 
$\H\ra \H/Z(\H)\simeq \R^{2n}$. 

We will begin by
introducing a quantity, $\alpha$, which
measures how close $E\in \fp(U)$ is to
being a half-space, taking into account location
and scale. For our purposes, not being close to a 
half-space is ``bad''.

Given a finite perimeter cut measure,
$\Si$,  we define the corresponding
the bad part,
$\lambda_{\eps,R}^{\bad}$,
of the total perimeter measure, $\lambda=\lambda_\Si$,
where the parameters, $\eps$, $R$ specify
the degree of badness and the scale respectively.
 Control on $\lambda_{\eps,R}^{\bad}$, which is
the key to proving
our main result,  is
obtained in Section \ref{badpercontrol}. 
Theorem \ref{maindiffcut},
is proved by translating the bounds 
on $\lambda_{\eps,R}^{\bad}$ into bounds
on the cut measure $\Si$.

\vskip4mm

\no
\subsection{Measuring closeness of FP cuts to half-spaces.}

\mbox{}

\vspace{2mm}We denote   the collection of all  half-spaces in $X$ by 
$\halfspace$,
and let
$$
\hs_x\defeq \{E\in \hs\mid \D E\;\mbox{contains}\;x\}\, .
$$

\begin{definition}
\label{alpha}
Define 
$$
\al:\fp(U)\times U\times (0,\infty)\ra \R_+
$$
to be the normalized $L^1$-distance between 
$E$ and $\hs_x$
in the ball $B_r(x)$:
$$
\al(E,x,r)= \min_{H\in \hs_x} \frac{1}{\mu(B_r(x))}
\int_{B_r(x)}|\chi_E-\chi_H|\, d\mu\, .
$$
\end{definition}
\vskip1mm

\begin{lemma}
$\al$ is a locally Lipschitz function of all $3$ variables.
\end{lemma}
\proof
Changing $r$ and $x$ slightly only adds or 
subtracts a small amount of measure.
Locally Lipschitz dependence on $E$ is clear.
\qed
\vskip2mm

 For $\eps,R>0$, and $E\in\fp(U)$, let

$$
\begin{aligned}
\bad_{\eps,R}(E)
&\defeq
\{x\in U\mid d(x,X\setminus U)<R\quad\mbox{or}\quad
\al(E,x,r)>\eps\;\mbox{for some}\;r\in(0,R]\}\, ,\\
&{}\\
\good_{\eps,R}(E)
&\defeq\{x\in U\mid d(x,X\setminus U)\leq R\quad\mbox{and}\quad\al(E,x,r)\leq\eps\;
\mbox{for all}\;r\in(0,R]\}\, .
\end{aligned}
$$

Thus,
$$
\good_{\eps,R}(E)=U\setminus\bad_{\eps,R}(E)\, .
$$
Also, put
$$
\begin{aligned}
\bad_{\eps,R}
&\defeq \{(E,x)\in \fp(U)\times U\mid
x\in \bad_{\eps,R}(E)\}\, ,\\
&{}\\
\good_{\eps,R}&\defeq \{(E,x)\in \fp(U)\times U\mid
x\in \good_{\eps,R}(E)\}\, . 
\end{aligned}
$$

\begin{lemma}
\label{open}
  $\bad_{\eps,R}$ is an open subset of $\fp(U)\times U$.  
\end{lemma}
\proof
The set $\bad_{\eps,R}$ is the image of the open set

$$
\{(E,x,r)\subset \fp(U)\times U\times (0,R]\mid 
d(x,X\setminus U)<R\quad\mbox{or}\quad \al(E,x,r)>\eps\}
$$
under the open projection map,
$$
  \fp(U)\times U\times (0,R]\ra \fp(U)\times U\, .
$$
The conclusion follows.
\qed
\vskip2mm

\no
\subsection{The total bad perimeter measure.}

\mbox{}

\vspace{2mm}As in Section \ref{totper}, let $\Si$ denote an FP cut measure
and $\lambda=\lambda_\Si$ the associated
total perimeter measure.

 Given a measure $\zeta$ on $Z$,
and  a  measurable subset $A\subset Z$, 
 let $\zeta\on A$ denote the measure given by
\begin{equation}
\label{measrestdef} 
\zeta\on A(F)=\zeta(A\cap F)\, .
\end{equation}

Let the map, 
$$
\per\on\bad_{\eps,R}:\fp(U)\lra \radon(U)\, , 
$$
be given by
\begin{equation}
\label{badpermeasdef}
\per\on\bad_{\eps,R}(E)=\per(E)\on\bad_{\eps,R}(E)\, .
\end{equation}

Recall the notion of weakly $L^1$ map; see (\ref{Psiellone}).
\begin{lemma}
The map 
$$
 \per\on\bad_{\eps,R}:\fp(U)\lra \radon(U)\,,
$$
is weakly $L^1$.
\end{lemma}
\proof
 For all $k$, let 
$$
\Phi_k:\fp(U)\times U\ra \R
$$
denote a continuous function satisfying:

1. $0\leq \Phi_k\leq 1$.

2. $\Phi_k\equiv 1$ on 
$$
\left\{(E,x)\mid 
d\left((E,x),\good_{\eps,R}\right)\geq \frac{1}{k}\right\}
\, ,
$$
where the distance on the product is  the sum 
of the factor distances.

3. $\Phi_k\equiv 0$ on 
$$
\left\{(E,x)\mid 
d\left((E,x),\good_{\eps,R}\right)\leq \frac{1}{k+1}\right\}
\,.
$$

 Fix $\phi\in C_c(U)$ and define
$$
\Psi_k:\fp(U)\ra \R
$$
 by 
$$
\Psi_k(E)=
\int_U \phi\;\Phi_k(E,\,\cdot\,)\;d\per(E)\, .
$$
The map $\Psi_k$ is Borel measurable,
since it is the pointwise limit of a sequence
of measurable functions  $\{\Psi_{k,l}\}$ 
obtained by approximating the map
$$
 E\to \phi(\,\cdot\,)\,\Phi_k(E,\,\cdot\,)
$$
by  simple functions, and 
each of the $\Psi_{k,l}$'s is measurable.

 For fixed $E$,  the compact subsets,
$$
\supp(\Phi_k(E,\,\cdot\,))\subset \bad_{\eps,R}(E),
$$
exhaust the open set $\bad_{\eps,R}(E)$. 
To see this note that
for each compact set, 
$K\subset \bad_{\eps,R}(E)$, the subset,
$$
\{E\}\times K\subset\fp\times U\, ,
$$
has positive distance from the closed 
set $\good_{\eps,R}(E)$, and is therefore
contained in $\supp(\Phi_k(E,\,\cdot\,))$ for 
sufficiently large $k$.

It follows from the above that the mass 
of the difference measure
$$
\per(E)\on\bad_{\eps,R}(E)-\Phi_k(E,x)\per(E)
$$
tends to $0$  as $k\ra \infty$.  
Thus, for each
  $E\in \fp(U)$, the integrals,
$$
\Psi_k(E)=\int_U \phi\;\Phi_k(E,x)\;d\per(E)\, ,
$$
converge as $k\ra\infty$, to 
$$
\int_U\phi\; d\left(\per(E)\on\bad_{\eps,R}(E)\right)\, .
$$

The map
$$
E\to \int_ U\phi \;
d\left(\per(E)\on\bad_{\eps,R}(E) \right)
$$
is a pointwise
limit of Borel measurable functions and 
is therefore Borel measurable. Since $\phi$ is arbitrary,
it follows that the map,
$$\per\on\bad_{\eps,R}:\fp(U)\to\radon(U)
$$
is weakly measurable.

Now (\ref{FPcutmeas}) implies that $\per\on\bad_{\eps,R}$
 is weakly $L^1$.

\qed

\begin{definition}
\label{deftotalbad}
The {\it total bad perimeter measure} is the Radon measure 
$$
\lambda_{\eps,R}^{\bad}\in\radon(U)
$$ 
obtained by applying 
Proposition \ref{lemtotalmeasure} to the weakly $L^1$
map 
$$
\per\on\bad_{\eps,R}:\fp(X)\to\radon(U)\,.
$$
\end{definition}

\section{Controlling the total bad perimeter measure}
\label{badpercontrol}

Recall that from now on $(X,\mu)$ will denote either $\R^n$ or $\H$, and $U\subset X$
will denote a ball.  Also, $\Si$ will denote an FP cut measure on $\cut(U)$, with associated
total perimeter measure $\la$,  and associated good and bad measures
$\la^{\good}_{\eps,R}$, $\la^{\bad}_{\eps,R}$.

One of the main results of \cite{italians1} 
(see \cite{degiorgi} for the $\R^n$ case) states that if $E$ 
is a set of finite perimeter
in $U$, then for $\per(E)$ a.e. $x\in U$, the
blow ups of $E$ at $x$ 
converge in  $L^1_{\loc}$ to a half-space.
In this section we give a version of the above result
in the parametrized setting. Namely, given
an FP cut measure $\Si$ with total perimeter measure,
$\lambda=\lambda_\Si$, we show
in Theorem \ref{badra0p}, that
for any fixed $\eps$, the mass of 
$\lambda_{\eps,R}^{\bad}$, goes to zero as
$R\to 0$.
Theorem \ref{badra0p} is of crucial importance;
its proof constitutes the one and only 
place where we explicitly appeal to \cite{italians1}.

 From Theorem \ref{badra0p} and a straightforward
differentiation argument, it follows that
 for any $\eps>0$
 we can find a set with almost full
measure on which 
$$
\frac{\lambda_{\eps,R}^{\bad}(B_r(x))}{\mu(B_r(x))}
$$
is as small as we like,
provided we take $R$ sufficiently small; see Proposition
\ref{xdeeps1}.

\begin{theorem}
\label{badra0p}
 For all $\eps>0$,
\begin{equation}
\label{badra0}
\lim_{R\to 0}\mass(\lambda_{\eps,R}^{\bad})= 0\, .
\end{equation}
\end{theorem}
\proof
By (\ref{eqntotalmass}), (\ref{badpermeasdef}),
\begin{equation}
\label{eqnmassbad}
\mass(\lambda_{\eps,R}^{\bad})
=\int_{\cut(U)}\mass\left(\per(E)\on(\bad_{\eps,R}(E)
\right)\, d\Si\, .
\end{equation}

By the main result of \cite{italians1},
 for fixed $\eps>0$, $E\in\fp$, we have
\begin{equation}
\label{mainital}
\lim_{R\ra 0}\mass(\per(E)\on(\bad_{\eps,R}(E))= 0\, .
\end{equation}
(For equivalent ways of expressing this, compare 
(\ref{massdef}) and (\ref{measrestdef}).)

 Since 
$$
\mass(\per(E)\on(\bad_{\eps,R}(E))
\leq \mass(\per(E))\, ,
$$
and $\per\in L^1(\cut(U),\Si)$ (see 
Definition \ref{deffiniteperimetercutmeasure} and 
Proposition \ref{bvcutchar}) the claim
 follows from (\ref{eqnmassbad}) and the dominated 
convergence theorem.
\qed

\bigskip
\begin{proposition}
\label{xdeeps1}
For all $\de>0$, $\eps>0$,
 there exists
$r_0(\de,\Si)>0$, $r_1(\de,\eps,\Si)>0$,
 $R_0=R_0(\de,\eps,\Si)>0$, 
  and a subset, $U_{\de,\eps}(\Si)\subset X$,
such that 
\begin{equation}
\mu(U\setminus U_{\de,\eps})<2\de
                          \left(1+\mass(\lambda)\right)\, ,
\end{equation}
\begin{equation}
\label{taucontrolled}
\frac{\lambda(B_r(x))}{\mu(B_r(x))}<\de^{-1}\, ,
\qquad{\rm if}\,\, x\in  U_{\de,\eps},\,\, r\leq r_0(\de,\Si)\, ,
\end{equation}
\begin{equation}
\label{taubadsmall}
\frac{\lambda_{\eps,R_0}^{\bad}(B_r(x))}{\mu(B_r(x))}<\eps\, ,
\qquad{\rm if}\,\,x\in  U_{\de,\eps}\,,\,\, r\leq r_1(\de,\eps,\Si)\, . 
\end{equation}
\end{proposition}
\proof
Given Theorem \ref{badra0p},
this is a straightforward application 
of measure differentiation.

By the Lebesgue decomposition theorem, there 
exists, $U'\subset U$, with $\mu(U\setminus U')=0$,
such that $\lambda$ is absolutely continuous
with respect to $\mu$ on $U'$.

Since
$$
\int_{U'}\frac{d\lambda}{d\mu}\;d\mu\leq \mass(\lambda)<\infty\, ,
$$
there exists  $U_1\subset U$ 
such that 
$$
\mu(U\setminus U_1)<2\de\mass(\lambda)\, ,
$$
$$
\frac{d\lambda}{d\mu}
<\frac{\de^{-1}}{2}\qquad\mbox{on}\,\, U_1\, .
$$

By measure differentiation, for $\mu$ a.e. $x\in U_1$, 
$$
\lim_{r\ra\infty}\frac{\lambda(B_r(x))}{\mu(B_r(x))}
=\frac{d\lambda}{d\mu}(x)\, .
$$
Therefore, there exists $r_0(\de,\Si)>0$ and 
$U_2\subseteq U_1$, such that for all
$0<r\leq r_1$, $x\in U_2$,
$$
\mu(U_1\setminus U_2)<\de\, ,
$$
and
$$
\frac{\lambda(B_r(x))}{\mu(B_r(x))}<\de^{-1}\, ,
\qquad{\rm if}\,\, x\in U_2,\,\, r\leq r_0(\de,\Si)\, .
$$

Since by (\ref{badra0}),
%  \begin{equation}
$$
\lim_{R\ra 0}\mass(\lambda_{\eps,R}^{\bad})=0\, ,
$$
%  \end{equation}
there exists $U_3\subset U$,
with 
$$
\mu(U\setminus U_3)<\frac{\de}{2}\, ,
$$
and
 $R_0(\de,\eps,\Si)>0$, such that

$$
\frac{d\lambda_{\eps,R_0}^{\bad}}{d\mu}<\frac{\eps}{2}\, 
\qquad{\rm on}\,\,U_3\, .
$$

As above, by using measure differentiation,
there exists $U_4\subset U_3$
and $r_1(\de,\eps,\Si)>0$, such that
%  \begin{equation}
$$
\mu(U_3\setminus U_4)<\frac{\de}{2}\, ,
$$
%  \end{equation}
and
%  \begin{equation}
$$
\frac{\lambda_{\eps,R_0}^{\bad}(B_r(x))}{\mu(B_r(x))}<\eps\, ,
\qquad{\rm if}\,\, x\in U_4,\,\, r\leq r_1(\de,\eps,\Si)\, .
$$
%  \end{equation}

Now take
$U_{\de,\eps}\defeq U_2\cap U_4$.
\qed
\bigskip

\section{Collections of good and bad cuts}
\label{gbcuts}

In this section, given an FP cut measure, $\Si$,
we introduce sets of good and bad cuts,
$\G$, $\B$, where as usual, we
take into account location and scale.
Estimates on $\G$ and $\B$, are derived
from   Proposition \ref{xdeeps1}.

In Section \ref{goodcutmeas}, using the 
set $\G$, we will construct a measure,
$\widehat{\Si}$, which is supported on 
half-spaces. In Section \ref{mainthm}, our 
main theorem is established by proving
that, for $\mu$ a.e. $x\in U$,
 in the limit as $r\to 0$, the normalized
$L^1$-distance between, $d_\Si$ and 
$d_{\widehat{\Si}}$ converges to zero.

 For $\de>0$, $\eps>0$, let $r_0(\de,\Si)$,
$r_1(\de,\eps,\Si)$, $R_0(\de,\eps,\Si)$, 
$U_{\de,\eps}$,
be as in Proposition \ref{xdeeps1}.

For all $x\in U$, $r>0$, we define
$$
\G(x,\de,\eps,r,\Si)\subset\fp(U)\, ,
$$
$$
{}
$$
 $$
% \B=
\B(x,\de,\eps,r,\Si)\subset\fp(U)\, ,
$$
by
\vskip1mm

\begin{equation}
\label{GB}
\begin{aligned}
\G(x,\de,\eps,r,\Si)&=\{E\in \fp(U)\,\mid\, 
\overline{B_r(x)}\cap \good_{\eps,R_0}(E)\neq\emptyset\}\, ,\\
&{}\\
% \end{equation}
%Let $B\subset\fp(X)$ denote the complementary set, 
\B(x,\de,\eps,r,\Si)&=\fp(U)\setminus \G(x,\de,\eps,r,\Si)\, .
\end{aligned}
\end{equation}

Note that 
$\B$ is an open 
subset of $\fp(U)$. In particular, $\G$ and $\B$ are both 
Borel sets.

\begin{proposition}
\label{bpersmall}
 Pick $\de>0$, $\eps>0$.
If $x\in U_{\de,\eps}$ and 
\begin{equation}
\label{rsmall1}
r<\min\left(\frac{R_0}{2},r_1\right)\, ,
\end{equation}
then the total perimeter of 
$\B(x,\de,\eps,r,\Si)\subset \fp(U)$ in $B_r(x)$ is bounded by 
$\de\cdot\mu(B_r(x))$:
\begin{equation}
\label{perbsmall}
\frac{1}{\mu(B_r(x))}\int_\B \per(E)(B_r(x))\, d\Si< \eps\, .
\end{equation}
\end{proposition}
\proof
The definition of $\B$ together with (\ref{rsmall1})
 implies that if $E\in \B$, 
then $\overline{B_r(x)}\subset \bad_{\eps,R_0}(E)$.
Hence,
\begin{equation}
\per(E)(B_r(x))
=\left(\per(E)\on\bad_{\eps,R_0}(E)\right)(B_r(x)).
\end{equation}
Therefore,
$$
\begin{aligned}
\int_\B\per(E)(B_r(x))\, d\Si
&=\int_\B(\per(E)\on \bad_{\eps,R_0}(E))\, d\Si\\
&\leq\int_{\cut(U)}
     \left(\per(E)\on\bad_{\eps,R_0}(E)\right)\, d\Si\\
&{}\\
&=\lambda_{\eps,R_0}^{\bad}(B_r(x))\\
&{}\\
&< \de\cdot\mu(B_r(x))\, ,
\end{aligned}
$$
where the last inequality follows from (\ref{taubadsmall}).
(Actually, we only used $x\in U_4\supset U_{\de,\eps}$,
for $U_4$ as in the proof of Proposition \ref{xdeeps1}.)
\qed
\vskip1mm

 For the next lemma, we need
a standard fact concerning $\H$. Namely, there
exists $c>0$  such that if 
$H\in\hs_x'$, for some $x'\in B_r(x)$, then
\begin{equation}
\label{hsper}
 c \cdot r^{-1}\mu(B_{2r})\leq\per(H)(B_{r}(x))\, .
\end{equation}
 This is 	easy to see, for example, by employing the
coarea formula.

\begin{proposition}
\label{hspermbound}
There  are constants 
 $\eps_0>0,\,c_0<\infty$, such that
if 
\begin{equation}
\eps<\eps_0\, ,
\end{equation}
\begin{equation}
\label{rsmall0}
r<\min\left(\frac{r_0}{2},\frac{R_0}{2}\right)\, ,
\end{equation}
then,
\begin{equation}
\label{sigcontrolled}
\Si(\G)\leq c_0 r\de^{-1}\, . 
\end{equation}
\end{proposition}
\proof
By the definition of $E\in \G$, there exists 
$x'\in \ol{B_r(x)}\cap \good_{\eps,R_0}$.  
By (\ref{rsmall0}), this implies that
for $\alpha$ as in Definition \ref{alpha},
we have
$$
\al(E,x',2r)\leq \eps\, .
$$
 Therefore,  for some half-space, $H\in \hs_{x'}$,
$$
\frac{1}{\mu(B_{2r}(x'))}\int_{B_{2r}(x')}|\chi_H-\chi_E|
\, d\mu\leq \eps \, .
$$

Thus, by (\ref{hsper}) and the lower semicontinuity of perimeter with respect to 
$L^1$ convergence,  there exists
$c_1>0$ and $\eps_0>0$ 
such that for $\eps<\eps_0$, 
$$
c_1(n) r^{-1}\mu(B_{2r})\leq\per(E)(B_{2r}(x))\, .
$$
Therefore,
$$
\begin{aligned}
\lambda(B_{2r}(x))
&\geq \int_\G \per(E)(B_{2r}(x))\, d\Si\\
& >c_1 r^{-1}\mu(B_{2r}(x)\Si(\G). 
\end{aligned}
$$
% \end{equation}
Hence,

$$
\begin{aligned}
\Si(\G)
&< \frac{\lambda(B_{2r}(x))}{c_1r^{-1}\mu(B_{2r}(x))}\\
&< \delta^{-1}\frac{\mu(B_{2r}(x))}{c_1r^{-1}\mu(B_{2r}(x))}\\
&=c_0\de^{-1} r\, ,
\end{aligned}
$$

where the second inequality is a consequence of 
(\ref{taucontrolled}).
\qed

\section{The approximating cut measure supported on half-spaces}
\label{goodcutmeas}

We retain the notation from the preceding section: $\Si$ denotes
an FP cut measure on a ball $U\subset X$, $x\in U$,  
$\de>0$, $\eps>0$, $r>0$, and $\G=\G(x,\de,\eps,r,\Si)$, $\B=\B(x,\de,\eps,r,\Si)$
denote the corresponding sets of good and bad cuts.

We now construct a cut measure, $\widehat\Si(x,\de,\eps,r,\Si) $, 
supported on the collection of half-spaces $\hs$. The measure, 
$\widehat\Si$,  will be constructed 
by ``straightening'' the  cuts in $\G$.
If $x\in U_{\de,\eps}$ and $r$
satisfies the smallness  conditions 
(\ref{rsmall1}), (\ref{rsmall0})
in Propositions \ref{bpersmall}, \ref{hspermbound},
then the  cut metric $d_{\widehat{\Si}}$
will be  close to $d_{\Si}$ 
in the normalized $L^1$ metric
 on $B_r(x)$.

\bigskip
\begin{lemma}
\label{bormap}
There is a  Borel map, 

$$
\ga:\G\ra \hs\subset\fp(U)\, ,
$$
sending each $E\in\G$ to a half-space, $\ga(E)\in\hs(\H)$,
 such that for some 
$x'\in \overline{B_r(x)}$, 
\begin{equation}
\label{close2r}
\frac{1}{\mu(B_{2r}(x'))}\int_{B_{2r}(x')}
|\chi_E-\chi_{\ga(E)}|\, d\mu<2\eps\, .
\end{equation}
\end{lemma}
\proof
Let

$$
W\defeq 
\{(E,x',H)\in \fp(U)\times \overline{B_r(x)}\times\hs\mid H\in\hs_{x'}\}\, .
$$
% \end{equation}
The collection of elements, $(E,x',H)\in W$, which satisfy 
\begin{equation}
\frac{1}{\mu(B_{2r}(x'))}\int_{B_{2r}(x')}|\chi_E-\chi_{H}|\, d\mu
<2\eps\, ,
\end{equation}
is open in $W$, and as a consequence of its
definition, maps to $\G$ under the projection, 
$\fp(U)\times \overline{B_r(x)}\times\hs\ra \fp(U)$.
Therefore, we can construct a Borel section 
of this projection over $\G$ 
as follows.  

Each $E\in \G$ lies in an open set, 
$U_E\subset \G$, over which one
has a section $\si_E:U_E\ra W$ 
whose $\hs$ component is constant.  
A countable collection
of these will cover $\G$. Hence, 
there is a countable disjoint cover 
by  Borel sets, $\{V_i\}$, such that 
each $V_i$ lies in a $U_E$.
We define $\si:\G\ra W$ by declaring that 
its restriction to $V_i$
agrees with  $\si_E$ resricted to $V_i$. 
\qed
\bigskip

We define the  Borel measure, $\widehat\Si(x,\de,\eps,r,\Si)$, 
to be the pushforward under $\ga$ of the measure
$\Si\on \G$, where $\G=\G(x,\de,\eps,r,\Si)$:
\begin{equation}
\label{Sihat}
\widehat\Si(x,\de,\eps,r,\Si)\defeq \ga_*(\Si\on \G)\, .
\end{equation}
It follows immediately
 that $\widehat\Si$ is supported on
$\hs\subset\fp(U)$.  Since  each $E\in \hs$ contributes uniformly
bounded measure and uniformly bounded perimeter in $U$, it follows from 
(\ref{sigcontrolled}) that $\hat\Si$ is an $\fp$ cut measure.

\bigskip
\section{Proof of the main theorem}
\label{mainthm}
Here we prove Theorem \ref{maindiffcut},
the main differentiation assertion for $\fp$ cut measures.
Theorem \ref{maindiffcut} immediately implies
Theorem \ref{thmmainfirst}.
 For convenience, we will assume $U$ is a ball in $\H$.
The argument also applies, mutatis mutandis, if $\H$ is
replaced by $\R^k$.

We retain the notation from  Section \ref{goodcutmeas}.

 For all $r> 0$, let $S_r:\H\ra \H$ denote 
an automorphism which scales by
$r$, and let $S_{x,r}:\H\ra \H$ be the composition 

$$
\H\stackrel{S_r}{\lra}\H\stackrel{l_x}{\lra}\H,
$$
where $l_x:\H\lra\H$ denotes left translation by $x\in \H$.
The pullback of a distance, $d$, under $S_{x,r}$
is denoted $S_{x,r}^*(d)$.

 As in (\ref{eqncutmetric}), let $d_\Si$  denote the
distance on $U$ associated to a cut measure $\Si$.

The  $L^1$  distance between metrics, $d,d'$, on $A\subset U$
is
\begin{equation}
\label{ellonedist}
\|d-d'\|_{L^1}\defeq
\int_{A\times A}|d(x_1,x_2)-d'(x_1,x_2)|\,d\mu\,d\mu\, .
\end{equation}

We let $\meas_{\hs}$ 
denote the collection of 
FP cut measures which are supported on half-spaces.

\begin{theorem}
\label{maindiffcut}
 Given an FP cut measure, $\Si$,
there is a subset 
$U_0\subset U$ of full Lebesgue measure 
such that if $x\in U_0$, then
 \begin{equation}
\label{convtohs}
\lim_{r\ra 0}\,\, \inf_{\bar\Si\in\meas_{\hs}}
\,\,\| \frac{1}{r}S_{x,r}^*(d_\Si)-d_{\bar\Si}\|_{L^1}= 0\, ,
{}
\end{equation}
where the $L^1$ norm is taken on the unit ball $B_1(e)$.
In particular, if $\Si=\Si_f$ is the cut measure corresponding
to a BV map, $f:\H\to L^1(Y)$, then (\ref{convtohs}) holds.
\end{theorem}
\proof 
Let $x,r$ satisfy the hypotheses of Propositions
\ref{bpersmall}, \ref{hspermbound} and Lemma \ref{bormap}.

We will show 
that on $B_r(x)$,
\begin{equation}
\label{final}
\|d_{\Si}-d_{\widehat\Si}\|_{L^1}
\leq 
r(4 c_0\eps\de^{-1}+\tau' \de )(\mu(B_r(x)))^2\, .
\end{equation}

Here $\tau'$ denotes the constant in the Poincar\'e
inequalty (\ref{pi2}).   The theorem follows by letting  
$\eps\to 0$ (which requires $r\to0$),
and then $\de\to 0$.  

Let $\widehat{\Si}$ denote the FP cut measure supported on 
half-spaces defined in (\ref{Sihat}).

On $B_r(x)$, the triangle inequality gives
\begin{equation}
\label{triangle}
\begin{aligned}
\|d_{\Si}-d_{\widehat\Si}\|_{L^1}
&\leq \|d_{\Si\on \G}-d_{\widehat\Si}\|_{L^1}+
\|d_{\Si\on \G}-d_{\Si}\|_{L^1}\\
&=\|d_{\Si\on \G}-d_{\widehat\Si}\|_{L^1}
+\|d_{\Si\on \B}\|_{L^1}\, .
\end{aligned}
\end{equation}
To complete the proof,
we estimate each term on the right-hand side on the second line of 
(\ref{triangle}).

The estimate from Proposition \ref{hspermbound}, which bounds the good cut measure, enters
in the proof of the next Lemma in a crucial way.   Without it we would only be able
to estimate the $L^1$ discrepancy between individual good cuts and their half-space
approximations, but would be unable to estimate the aggregate effect 
on the cut metric of this  discrepancy.

\begin{lemma}
\label{lemma1}
On $B_r(x)$,
\begin{equation}
\label{e}
\|d_{\Si\on \G}-d_{\widehat\Si}\|_{L^1}
\leq 4 c_0r\eps\de^{-1}\, .
\end{equation}
\end{lemma}
\proof
On $B_r(x)$, the  left-hand side  of (\ref{e}) is equal to 
$$
\begin{aligned}
&\int_{B(x,r)\times B_r(x)}|d_{\Si\on \G}(x_1,x_2)
-d_{\widehat\Si}(x_1,x_2)|\, d\mu\, d\mu\\
&\leq \int_{B_r(x)\times B_r(x)}
\int_\G\left||\chi_E(x_1)-\chi_E(x_2)|
-|\chi_{\ga(E)}(x_1)-\chi_{\ga(E)}(x_2)|\right|
\, d\Si\,d\mu\,  d\mu\\
&\leq \int_\G\int_{B_r(x)} 2\mu(B_r(x))|\chi_E-\chi_{\ga(E)}|\,d\mu
\, d\Si\\
&\leq 4\int_\G  \mu(B_r(x))(\eps\mu(B_r(x)))
\, d\Si\quad \mbox{(by (\ref{close2r}))} \\
&\leq 4 \eps\Si(\G)(\mu(B_r(x)))^2\\
&\leq 4 c_0r\eps\de^{-1}\quad\mbox{(by  (\ref{sigcontrolled}))}\,.
\end{aligned}
$$
\qed

\begin{lemma}
\label{lemma2}
On $B_r(x)$,
\begin{equation}
\label{db}
\|d_{\Si\on \B}\|_{L^1}
\leq r\tau \de (\mu(B_r(x)))^2\, .
\end{equation}
\end{lemma}
\proof
\begin{equation}
\label{ee}
\begin{aligned}
\|d_{\Si\on \B}\|_{L^1}
&=\int_B\int_{B_r(x)\times B_r(x)}
\, d_E(x,y)\, d\mu(x)\, d\mu(y)\, d\Si\\
&\leq cr\mu(B_r(x))\int_\B\per(E)(B_r(x))\,d\Si\,\, ,
\end{aligned}
\end{equation}
where the last inequality follows from the Poincar\'e inequality
(\ref{pi2}).
 From (\ref{ee}) and (\ref{perbsmall}), we get (\ref{db}).
\qed

Combining (\ref{triangle}),
(\ref{e}) and (\ref{db}) gives (\ref{final}),
which suffices to complete the proof.
\qed
\vskip1mm

\begin{remark}
 A refinement of the proof of Theorem \ref{maindiffcut}
yields the following stronger statement: 
 For $\mu$ a.e. $x\in U$,
blow ups of the FP cut measure converge
to a translation invariant cut measure 
which is supported on half-spaces.
\end{remark}
\vskip1mm

\begin{remark}
  The proof presented
here works for any Carnot group $G$ for which the blow up
result of \cite{italians1} is valid; in particular, it holds
for an arbitrary 2-step nilpotent 
Lie group.   It seems  almost certain that
their result will hold for general nilpotent Lie groups.  
\end{remark}

\bibliography{ellone}
\bibliographystyle{alpha}

\end{document}